\documentclass[preprint,12pt,authoryear]{elsarticle}

\usepackage{amssymb}
\usepackage{amsmath} 
\usepackage{graphicx}
\usepackage{hyperref}
\usepackage{booktabs}
\usepackage{listings}
\usepackage{xcolor}
\usepackage{ragged2e}
\usepackage{float}
\usepackage{subcaption}
\usepackage{algorithm} 
\usepackage{algpseudocode} 
\usepackage{caption}
\usepackage{setspace}
\usepackage{amsmath,amssymb}
\usepackage{bm}
\usepackage{multirow}
\usepackage[utf8]{inputenc}
\usepackage[modulo]{lineno}
\usepackage{tabularx}
\usepackage{booktabs}
\usepackage{makecell}
\journal{}

\begin{document}

\begin{frontmatter}



\title{Nonlinear Geotechnical Analysis Using a Polygonal Cell-Based Smoothed Finite Element Framework} %% Article title

\author[inst2]{Mingjiao Yan}
\affiliation[inst2]{organization={College of Water Conservancy and Hydropower Engineering, Hohai University},%Department and Organization
            %addressline={Address One}, 
            city={Nanjing},
            postcode={210098}, 
            state={Jiangsu},
            country={China}}

\author[inst1]{Yang Yang}
\affiliation[inst1]{organization={PowerChina Kunming Engineering Corporation Limited},%Department and Organization
            %addressline={Address Two}, 
            city={Kunming},
            postcode={650051}, 
            state={Yunnan},
            country={China}}

\affiliation[inst3]{organization={College of Harbour, Coastal and Offshore Engineering, Hohai University},%Department and Organization
            %addressline={Address Two}, 
            city={Nanjing},
            postcode={210098}, 
            state={Jiangsu},
            country={China}}

\author[inst1]{Zongliang Zhang}
\author[inst3]{Yinpeng Yin}
\author[inst2]{Miao Zhang}
\author[inst2]{Yijia Dong}
\author[inst1]{Dong Pan}
\author[inst1]{Xiaozi Lin}
\author[inst1]{Tiankai Yang}
%% Abstract
\begin{abstract}
Nonlinear geotechnical analysis often involves complex geometries, staged construction, local failure, and mesh-dependent stress and plastic strain responses. This study develops a polygonal cell-based smoothed finite element method (CS-FEM) for nonlinear geotechnical analysis and implements it in ABAQUS through the user element subroutine. The proposed method combines Wachspress interpolation with cell-based strain smoothing, in which the smoothed strain--displacement matrix is evaluated by boundary integration over polygonal smoothing subcells. This formulation avoids direct calculation of shape-function derivatives inside polygonal elements and enables standard polygonal meshes and hybrid quadtree meshes with hanging nodes to be handled in a unified framework. Nonlinear geomaterial behavior is incorporated through incremental elasto-plastic constitutive updates, including the Mohr--Coulomb model and the Duncan--Chang model. Several benchmark and engineering examples, including a perforated plate, strip footing, core rockfill dam, tunnel excavation, and slope stability problems, are presented for verification. The results show that the proposed method accurately predicts displacement, stress, plastic strain, bearing capacity, and factor of safety, while providing improved mesh flexibility and computational efficiency for nonlinear geotechnical analysis.
\end{abstract}

\begin{highlights}
\item A polygonal CS-FEM is developed for nonlinear geotechnical analysis.
\item Wachspress interpolation is combined with cell-based strain smoothing.
\item Boundary integration avoids shape-function derivatives in polygonal cells.
\item ABAQUS UEL enables incremental elasto-plastic analysis of geomaterials.
\end{highlights}

%% Keywords

\end{frontmatter}



\section{Introduction}
\label{sec:Introduction}

The numerical prediction of geomaterial response under complex loading conditions remains a central topic in computational geomechanics. This difficulty arises from the intrinsically nonlinear and history-dependent behavior of soils and rocks, which may involve plastic yielding, pressure-dependent strength, strain softening, progressive failure, time-dependent deformation, and pronounced loading-path effects~\citep{castellon2022small,zhang2021thermodynamics,liu2023unified,wang2023progressive,li2025progressive}. Reliable modeling of these behaviors is essential for the analysis and design of a wide range of geotechnical systems, including foundations, slopes, tunnels, embankments, and underground structures. Owing to its generality and versatility, the finite element method (FEM) has become one of the most widely used computational tools for geotechnical analysis, particularly for problems involving complex geometries, material heterogeneity, nonlinear constitutive behavior, and nontrivial boundary conditions~\citep{zienkiewicz2005fem,bathe2006finite,hughes2000fem}.

Despite its widespread use, the conventional displacement-based FEM still faces several difficulties in nonlinear geotechnical analysis. Low-order finite elements may exhibit overly stiff responses under certain conditions, leading to underestimated displacement fields and reduced accuracy in stress prediction~\citep{fraeijs1965displacement,liu2003fem}. In addition, the use of isoparametric mapping requires valid Jacobians and imposes restrictions on element distortion, which may reduce the robustness of the method for irregular meshes, locally refined discretizations, and complex engineering geometries~\citep{knupp2000mesh,gargallo2015distortion}. These issues become more pronounced in nonlinear simulations of geomaterials, where plastic deformation, stress redistribution, and progressive failure are strongly affected by the quality and flexibility of the spatial discretization.

The smoothed finite element method (S-FEM), originally proposed by Liu et al.~\citep{liu2007sfem}, provides an effective strategy for improving the performance of standard FEM formulations. By introducing a strain-smoothing operation, the compatible strain field is replaced by a smoothed strain field evaluated through boundary integration over smoothing domains. This treatment can alleviate the overly stiff behavior of conventional FEM and reduce the dependence on shape-function derivatives and isoparametric mapping~\citep{liu2010smoothed,zeng2018sfemreview,dai2007n}. Different SFEM variants have been developed, including node-based, edge-based, face-based, and cell-based formulations~\citep{liu2009edge,nguyenthoi2010node,nguyen2009face,yan2025fast}. Among them, the cell-based smoothed finite element method (CS-FEM) is particularly attractive for engineering computation because its smoothing domains are constructed locally within each element, making the formulation simple, robust, and naturally compatible with standard finite element assembly procedures~\citep{natarajan2015equivalence,zeng2016effective}.

In parallel, polygonal and polyhedral finite element methods have attracted increasing attention as flexible alternatives to conventional triangular, quadrilateral, tetrahedral, and hexahedral discretizations. Polygonal elements offer several advantages, including enhanced geometric flexibility, improved mesh quality, convenient treatment of hanging nodes, and natural compatibility with adaptive refinement and complex domain decomposition~\citep{sukumar2006recent,talischi2012polymesher}. Shape functions over arbitrary polygons can be constructed using generalized barycentric coordinates, such as Wachspress coordinates~\citep{wachspress1975rational}, mean value coordinates~\citep{floater2003mean}, and natural neighbor interpolation~\citep{sukumar2001natural}. Nevertheless, the use of polygonal elements in nonlinear geotechnical analysis remains challenging. The numerical integration of rational shape functions, the treatment of nonlinear constitutive updates, and the stability of stress and plastic strain predictions require careful formulation. Although polygonal SFEM formulations have been investigated for linear problems~\citep{dai2007n,yang2026seepage}, their extension to nonlinear geomaterial analysis and engineering-scale geotechnical applications is still limited.

Another practical issue is that most S-FEM formulations are not directly available in commercial finite element software. Their application to engineering problems therefore often requires user-defined element implementation, for example through the ABAQUS UEL interface~\citep{abaqus2010user}. Compared with edge-based or node-based smoothing schemes, CS-FEM is more suitable for such an implementation because all smoothing operations are performed within each element and no additional inter-element topological reconstruction is required~\citep{kumbhar2020development,cui2020abaqus,colombo2023element}. This element-wise feature allows the proposed formulation to be incorporated into the nonlinear solution framework of ABAQUS while retaining the flexibility of polygonal discretization.

Motivated by these considerations, this study develops a polygonal cell-based smoothed finite element method for nonlinear geotechnical analysis and implements it in ABAQUS through the UEL interface. The proposed method combines Wachspress interpolation with cell-based strain smoothing, so that the smoothed strain--displacement matrix is evaluated by boundary integration over polygonal smoothing subcells. Nonlinear geomaterial behavior is incorporated through elasto-plastic constitutive updates, and the element tangent stiffness matrix and internal force vector are consistently formulated for incremental nonlinear analysis. The main contributions of this work are as follows: (1) a polygonal CS-FEM formulation is established for nonlinear elasto-plastic geotechnical analysis; (2) a local smoothing-domain construction is developed for polygonal elements using Wachspress interpolation and boundary integration; (3) the proposed formulation is implemented in ABAQUS via UEL, together with a post-processing procedure for stress and strain recovery; and (4) the accuracy, robustness, and engineering applicability of the method are verified through elastic, bearing capacity, staged construction, and slope stability problems.

The remainder of this paper is organized as follows. Section~\ref{sec:Cell-based SFEM formulation for elasto-plasticity} presents the polygonal cell-based smoothed finite element formulation, including Wachspress interpolation, strain smoothing, and the construction of the smoothed strain--displacement matrix. Section~\ref{sec:nonlinear} introduces the nonlinear elasto-plastic constitutive formulation used for geomaterial analysis. Section~\ref{sec:Implementation} describes the ABAQUS UEL implementation and the post-processing procedure for stress and strain recovery. Section~\ref{sec:numerical_example} presents several numerical examples, including an infinite plate with a circular hole, a strip footing bearing capacity problem, staged construction analysis of a core rockfill dam, and slope stability analysis. Finally, Section~\ref{sec:conclusions} summarizes the main conclusions and discusses future research directions.

\section{Polygonal cell-based smoothed finite element formulation}
\label{sec:Cell-based SFEM formulation for elasto-plasticity}

In this study, we examine a deformable body occupying the domain $\Omega$. The body is subjected to body forces $\mathbf{b}$, prescribed tractions $\overline{\mathbf{t}}$ on the boundary $\Gamma_t$, and displacement boundary conditions $\mathbf{u}=\overline{\mathbf{u}}$ on $\Gamma_u$. The governing equations are given by~\citep{zienkiewicz2005fem,bathe2006finite}
\begin{equation}
\begin{array}{rlc}
\nabla\cdot\boldsymbol{\sigma}+\mathbf{b}=0 & \text { in } & \Omega \\
\mathbf{u}=\overline{\mathbf{u}}  & \text{ on } & \Gamma_u, \\
\boldsymbol{\sigma}\mathbf{n}=\overline{\mathbf{t}} & \text{ on } & \Gamma_t
\end{array}
\end{equation}
where $\boldsymbol{\sigma}$ denotes the Cauchy stress tensor.

The constitutive relation is expressed as
\begin{equation}
\boldsymbol{\sigma}=\mathbf{D}\boldsymbol{\varepsilon},
\end{equation}
where $\mathbf{D}$ is the constitutive matrix.

The computational domain is discretized into arbitrary convex polygonal elements. Within each element, the displacement field is approximated using Wachspress shape functions defined over polygons as
\begin{equation}
\mathbf{u}^h(\mathbf{x})=\sum_{i=1}^{n} N_i(\mathbf{x}) \mathbf{d}_i,
\end{equation}
where $n$ is the number of polygon nodes and $N_i$ are Wachspress interpolants.

The Wachspress interpolants are adopted to approximate the displacement field over a convex polygonal element. They can be expressed in the normalized rational form as
\begin{equation}
N_i(\mathbf{x})=\frac{w_i(\mathbf{x})}{\sum_{j=1}^{n} w_j(\mathbf{x})},
\end{equation}
where $w_i(\mathbf{x})$ denotes the Wachspress weight associated with node $i$. These shape functions satisfy the partition of unity, Kronecker delta property, and linear completeness, and therefore are suitable for constructing conforming polygonal finite elements.

The virtual displacement and strain fields are given by
\begin{equation}
\delta \mathbf{u}^h=\sum_i N_i \delta \mathbf{d}_i, 
\quad 
\mathbf{u}^h=\sum_i N_i \mathbf{d}_i,
\end{equation}
\begin{equation}
\delta \boldsymbol{\varepsilon}^h=\sum_i \mathbf{B}_i \delta \mathbf{d}_i, 
\quad 
\boldsymbol{\varepsilon}^h=\sum_i \mathbf{B}_i \mathbf{d}_i,
\end{equation}
where $\mathbf{d}_i=[u_i\ v_i]^{\mathrm{T}}$, and the strain–displacement matrix in 2D is
\begin{equation}
\mathbf{B}_i=
\begin{bmatrix}
N_{i,x} & 0 \\
0 & N_{i,y} \\
N_{i,y} & N_{i,x}
\end{bmatrix}.
\end{equation}

Applying the principle of virtual work yields~\citep{liu2007theoretical}
\begin{equation}
\int_{\Omega} \mathbf{B}^{\mathrm{T}} \boldsymbol{\sigma} \, \mathrm{d}\Omega
=
\int_{\Omega} \mathbf{N}^{\mathrm{T}}\mathbf{b} \, \mathrm{d}\Omega
+
\int_{\Gamma_t} \mathbf{N}^{\mathrm{T}} \mathbf{t} \, \mathrm{d}\Gamma.
\end{equation}

The discrete system can be written as
\begin{equation}
\mathbf{K}\mathbf{d}=\mathbf{P},
\end{equation}
with
\begin{equation}
\mathbf{K}=\int_{\Omega} \mathbf{B}^{\mathrm{T}} \mathbf{D}\mathbf{B} \, \mathrm{d}\Omega,
\quad
\mathbf{P}=\int_{\Omega} \mathbf{N}^{\mathrm{T}} \mathbf{b} \, \mathrm{d}\Omega+\int_{\Gamma_t} \mathbf{N}^{\mathrm{T}} \mathbf{t} \, \mathrm{d}\Gamma.
\end{equation}

\subsection{Strain smoothing in polygonal CS-FEM}

In the CS-FEM framework, the strain field is smoothed over each smoothing sub-cell $\Omega_C$ as
\begin{equation}
\tilde{\boldsymbol{\varepsilon}}^h
=\int_{\Omega_C} \boldsymbol{\varepsilon}^h(\mathbf{x}) 
\Phi(\mathbf{x}-\mathbf{x}_C) \, \mathrm{d}\Omega,
\end{equation}
where the smoothing function is defined as
\begin{equation}
\Phi(\mathbf{x}-\mathbf{x}_C)=
\begin{cases}
1/|\Omega_C|, & \mathbf{x}\in\Omega_C \\
0, & \text{otherwise}
\end{cases}.
\end{equation}

Applying the divergence theorem leads to
\begin{equation}
\tilde{\boldsymbol{\varepsilon}}^h(\mathbf{x}_C)
=
\frac{1}{|\Omega_C|}
\int_{\Gamma_C} \mathbf{u}^h(\mathbf{x}) \mathbf{n}(\mathbf{x}) \, \mathrm{d}\Gamma.
\end{equation}

Substituting the displacement approximation gives
\begin{equation}
\tilde{\boldsymbol{\varepsilon}}^h
=
\sum_{i=1}^{n} \widetilde{\mathbf{B}}_i \mathbf{d}_i,
\end{equation}
where the smoothed strain–displacement matrix is
\begin{equation}
\widetilde{\mathbf{B}}_i
=
\frac{1}{|\Omega_C|}
\int_{\Gamma_C}
\begin{bmatrix}
n_x & 0 \\
0 & n_y \\
n_y & n_x
\end{bmatrix}
N_i(\mathbf{x}) \, \mathrm{d}\Gamma.
\end{equation}

Thus, the stiffness matrix in CS-FEM is obtained as
\begin{equation}
\widetilde{\mathbf{K}}
=
\sum_{e}
\int_{\Omega^e}
\widetilde{\mathbf{B}}^{\mathrm{T}}
\mathbf{D}
\widetilde{\mathbf{B}} \, \mathrm{d}\Omega.
\label{eq:kmatrix}
\end{equation}

\subsection{Supported element types in the polygonal CS-FEM formulation}
\label{sec:supported_element_types}

\begin{figure}[H] 
\centering 
\includegraphics[width=0.9\textwidth]{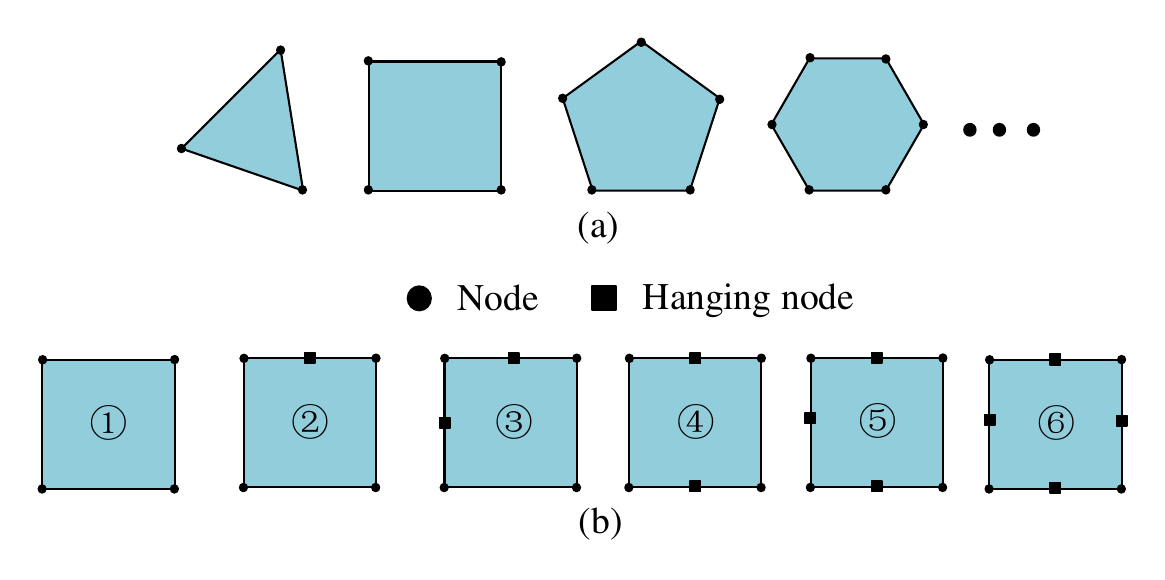} 
\caption{Supported element types in the proposed polygonal CS-FEM formulation: (a) polygonal elements with different numbers of vertices; (b) quadrilateral elements with different hanging-node arrangements.}
\label{fig:csfem_element_types}
\end{figure}

The proposed polygonal CS-FEM formulation can be applied to different types of two-dimensional elements. As shown in Fig.~\ref{fig:csfem_element_types}, the supported elements include standard polygonal elements and quadtree elements with hanging nodes. For standard polygonal elements, the number of vertices is not limited to three or four, and the displacement field can be interpolated over each convex polygon using Wachspress shape functions.

For locally refined quadtree meshes, hanging nodes may appear on element edges when adjacent elements have different refinement levels. In the present formulation, the hanging nodes are treated as regular boundary nodes and included in the ordered nodal list of the element. Therefore, quadrilateral elements with different numbers and arrangements of hanging nodes can be regarded as polygonal cells and handled in the same framework without introducing transition elements or additional constraint equations.

\subsection{Construction of smoothing subcells}
\label{sec:smoothing_subcells}

\begin{figure}[H] 
\centering 
\includegraphics[width=0.85\textwidth]{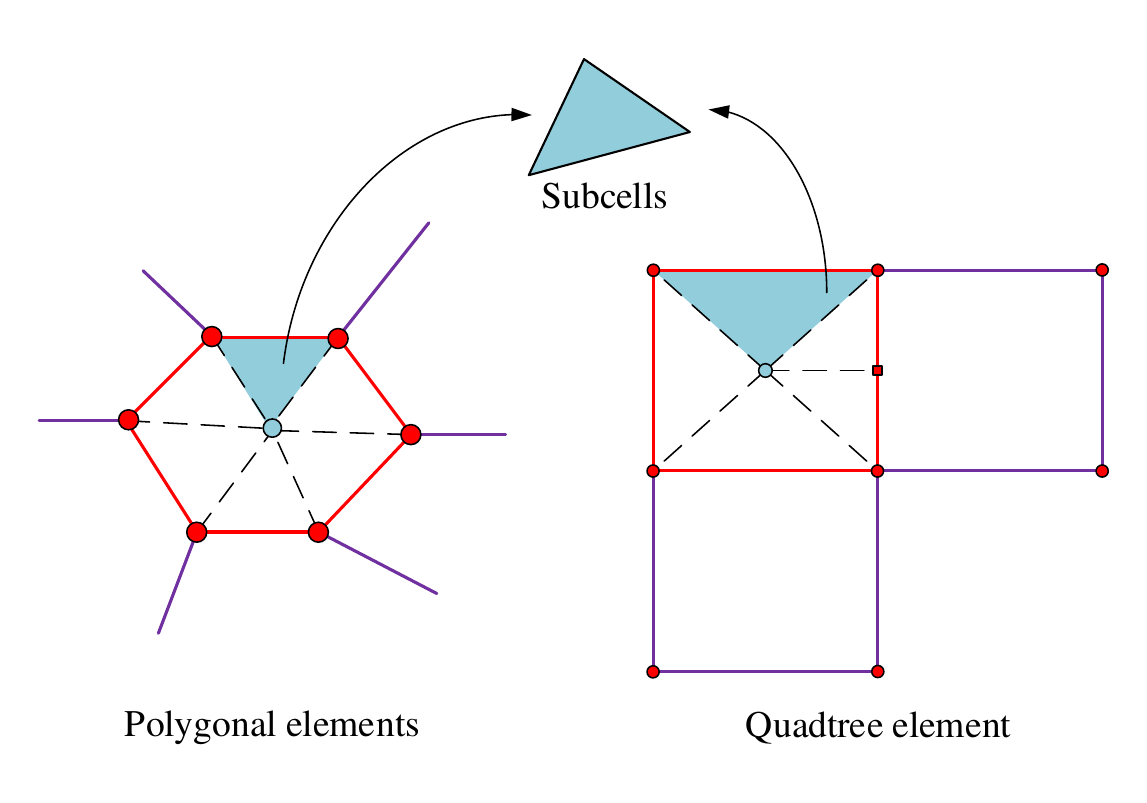} 
\caption{Construction of smoothing subcells for polygonal and quadtree elements. For both element types, the element centroid is connected to the boundary nodes to generate triangular smoothing subcells.}
\label{fig:csfem_ele_gen}
\end{figure}

In the present polygonal CS-FEM formulation, each element is further divided into a set of smoothing subcells for strain smoothing. As shown in Fig.~\ref{fig:csfem_ele_gen}, the same centroid-based subdivision strategy is adopted for both general polygonal elements and quadtree elements with hanging nodes. For a general polygonal element, the smoothing subcells are constructed by connecting the element centroid to the element vertices. Each subcell is bounded by one element edge and two line segments connecting the centroid to the two end nodes of that edge.

For quadtree elements with hanging nodes, the hanging nodes located on element edges are included in the element boundary-node sequence. The element centroid is then connected to all boundary nodes, including both corner nodes and hanging nodes, so that the element is subdivided into several triangular smoothing subcells. In this way, both standard polygonal elements and quadtree elements with hanging nodes can be processed using the same smoothing-cell construction procedure.

\subsection{Wachspress shape functions for polygonal elements}
\label{sec:wachspress_poly}

After the smoothing subcells are constructed, all supported elements, including standard polygonal elements and quadtree elements with hanging nodes, can be regarded as convex polygonal cells with ordered boundary nodes. Wachspress shape functions~\citep{wachspress1975rational} are then employed to evaluate the displacement interpolation and the boundary integrals over the smoothing subcell edges.

\begin{figure}[H] 
\centering 
\includegraphics[width=0.65\textwidth]{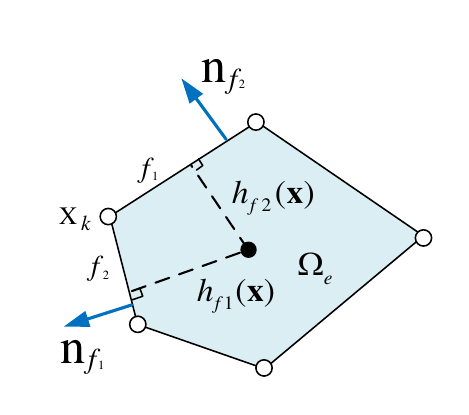} 
\caption{Geometric quantities used in the construction of Wachspress shape functions on a convex polygonal element.}
\label{fig:shape_function}
\end{figure}

Each polygonal element $\Omega_e$ is defined as a convex $n$-gon with ordered vertices $\{\mathbf{x}_k\}_{k=1}^{n}$ in the physical plane. Fig.~\ref{fig:shape_function} schematically illustrates the construction of Wachspress shape functions over such a polygon.

For a given vertex $k$, let $f_1$ and $f_2$ denote the two edges adjacent to it, and $\mathbf{n}_{f_i}$ be the outward unit normal vectors of these edges. The signed distance from any point $\mathbf{x}\in \Omega_e$ to the line supporting edge $f_i$ is defined as
\begin{equation}
h_{f_i}(\mathbf{x}) = (\mathbf{x} - \mathbf{x}_{f_i}) \cdot \mathbf{n}_{f_i},
\end{equation}
where $\mathbf{x}_{f_i}$ is a point on edge $f_i$ (e.g., one endpoint). By convention, $h_{f_i}(\mathbf{x})>0$ inside the polygon.

The Wachspress weight associated with vertex $k$ is then given by
\begin{equation}
w_k(\mathbf{x}) = \frac{\det(\mathbf{n}_{f_1},\mathbf{n}_{f_2})}{h_{f_1}(\mathbf{x})\, h_{f_2}(\mathbf{x})},
\end{equation}
where the determinant depends only on the local geometry of edges adjacent to vertex $k$.

Using these weights, the Wachspress shape function for vertex $k$ is constructed as
\begin{equation}
N_k(\mathbf{x}) = \frac{w_k(\mathbf{x})}{\sum_{j=1}^{n} w_j(\mathbf{x})}.
\label{eq:shape_function_values_2d}
\end{equation}

These shape functions are nonnegative, $C^0$-continuous over the polygon, and satisfy the linear completeness and partition-of-unity properties:
\begin{equation}
\sum_{k=1}^{n} N_k(\mathbf{x}) = 1 \quad \forall \mathbf{x} \in \Omega_e.
\end{equation}

In the CS-FEM framework, Wachspress functions are used only in the boundary integrals over the smoothing subcells. For a subcell $C$ associated with node $I$, the smoothed strain–displacement contribution can be expressed as
\begin{equation}
\widetilde{\mathbf{B}}_I^{(C)} = \frac{1}{A_C} \int_{\Gamma_C} N_I(\mathbf{x})\, \mathbf{n}(\mathbf{x}) \, \mathrm{d}\Gamma,
\end{equation}
or equivalently, as a sum over edges of the subcell:
\begin{equation}
\widetilde{\mathbf{B}}_I^{(C)} = \frac{1}{A_C} \sum_{e=1}^{m_C} \left( \int_e N_I(\mathbf{x}) \, \mathrm{d}\Gamma \right) \mathbf{n}_e,
\label{eq:bk_2d}
\end{equation}
where $m_C$ is the number of edges of the subcell, $\mathbf{n}_e$ is the outward normal of edge $e$, and $A_C$ is the subcell area. Since the shape function $N_I$ is evaluated only along straight edges, the integrand is smooth and can be accurately integrated with a small number of Gauss points.

\section{Nonlinear constitutive formulation in polygonal CS-FEM}
\label{sec:nonlinear}

In nonlinear geomaterial analysis, the constitutive response may be described by either elasto-plastic models or nonlinear elastic models, depending on the material behavior of interest. In the present polygonal CS-FEM framework, the constitutive update is performed at the material points of each smoothing subcell, and the resulting stress and tangent constitutive matrix are used to assemble the element internal force vector and tangent stiffness matrix. In this study, the Mohr--Coulomb model is used for elasto-plastic soil behavior, while the Duncan--Chang hyperbolic model is used to describe the nonlinear stress-dependent stiffness of rockfill materials.

For the elasto-plastic formulation, the total strain increment at each material point is decomposed into elastic and plastic components as
\begin{equation}
\Delta \boldsymbol{\varepsilon} =
\Delta \boldsymbol{\varepsilon}^{\mathrm{e}}
+
\Delta \boldsymbol{\varepsilon}^{\mathrm{p}},
\end{equation}
where $\Delta \boldsymbol{\varepsilon}^{\mathrm{e}}$ and $\Delta \boldsymbol{\varepsilon}^{\mathrm{p}}$ denote the elastic and plastic strain increments, respectively. The elastic response follows Hooke's law:
\begin{equation}
\Delta \boldsymbol{\sigma}
=
\mathbf{D}\Delta \boldsymbol{\varepsilon}^{\mathrm{e}}
=
\mathbf{D}
\left(
\Delta \boldsymbol{\varepsilon}
-
\Delta \boldsymbol{\varepsilon}^{\mathrm{p}}
\right),
\end{equation}
where $\mathbf{D}$ is the elastic constitutive matrix.

The plastic strain increment is governed by the flow rule
\begin{equation}
\Delta \boldsymbol{\varepsilon}^{\mathrm{p}}
=
\frac{\partial g}{\partial \boldsymbol{\sigma}}
\Delta \lambda,
\end{equation}
where $g(\boldsymbol{\sigma},\kappa)$ is the plastic potential function, $\Delta \lambda$ is the plastic multiplier, and $\kappa$ denotes the internal hardening variable. During plastic loading, the stress state satisfies the consistency condition
\begin{equation}
\Delta F
=
\frac{\partial F}{\partial \boldsymbol{\sigma}}
\Delta \boldsymbol{\sigma}
+
\frac{\partial F}{\partial \kappa}
\Delta \kappa
=
0,
\end{equation}
where $F(\boldsymbol{\sigma},\kappa)$ is the yield function.

For a general non-associated elasto-plastic model, the elasto-plastic tangent constitutive matrix can be written as
\citep{sloan1987substepping}
\begin{equation}
\mathbf{D}_{ep}
=
\mathbf{D}
-
\frac{
\mathbf{D}\mathbf{a}_g \mathbf{a}_f^{\mathrm{T}}\mathbf{D}
}{
\mathbf{a}_f^{\mathrm{T}}\mathbf{D}\mathbf{a}_g + H
},
\label{eq:CSFEM_Dep}
\end{equation}
where
\begin{equation}
\mathbf{a}_f =
\frac{\partial F}{\partial \boldsymbol{\sigma}},
\quad
\mathbf{a}_g =
\frac{\partial g}{\partial \boldsymbol{\sigma}},
\end{equation}
and $H$ is the hardening modulus. For perfectly plastic materials, $H=0$. In the case of associated plastic flow, $g=F$ and thus $\mathbf{a}_g=\mathbf{a}_f$. This formulation is used for the Mohr--Coulomb model in the present study.

For the Duncan--Chang hyperbolic model, the material response is described by a stress-dependent nonlinear elastic tangent matrix. Therefore, no plastic multiplier or yield correction is required. The stress increment is written as
\begin{equation}
\Delta \boldsymbol{\sigma}
=
\mathbf{D}_{DC}
\Delta \boldsymbol{\varepsilon},
\end{equation}
where $\mathbf{D}_{DC}$ denotes the tangent constitutive matrix determined from the current stress state and the Duncan--Chang model parameters. In this way, both the Mohr--Coulomb elasto-plastic model and the Duncan--Chang nonlinear elastic model can be incorporated into the same polygonal CS-FEM solution framework by using an appropriate tangent constitutive matrix.

Within each polygonal element, the smoothed strain--displacement matrix $\widetilde{\mathbf{B}}$ is used to compute the element internal force vector:
\begin{equation}
\mathbf{F}_{\mathrm{int}}^e
=
\int_{\Omega^e}
\widetilde{\mathbf{B}}^{\mathrm{T}}
\boldsymbol{\sigma}
\,\mathrm{d}\Omega .
\end{equation}
The corresponding element tangent stiffness matrix is given by
\begin{equation}
\mathbf{K}_t^e
=
\int_{\Omega^e}
\widetilde{\mathbf{B}}^{\mathrm{T}}
\mathbf{D}_t
\widetilde{\mathbf{B}}
\,\mathrm{d}\Omega ,
\end{equation}
where $\mathbf{D}_t$ denotes the tangent constitutive matrix. Specifically, $\mathbf{D}_t=\mathbf{D}_{ep}$ for the Mohr--Coulomb elasto-plastic model, and $\mathbf{D}_t=\mathbf{D}_{DC}$ for the Duncan--Chang hyperbolic model.

The smoothed strain is obtained through boundary integration over the smoothing subcell $\Omega_C$ as
\begin{equation}
\tilde{\boldsymbol{\varepsilon}}^h
=
\frac{1}{|\Omega_C|}
\int_{\Gamma_C}
\mathbf{u}^h(\mathbf{x})
\mathbf{n}(\mathbf{x})
\,\mathrm{d}\Gamma .
\end{equation}

The global nonlinear equilibrium equation is assembled as
\begin{equation}
\mathbf{R}
=
\mathbf{F}_{\mathrm{ext}}
-
\sum_e
\mathbf{F}_{\mathrm{int}}^e
=
\mathbf{0}.
\end{equation}

The system is solved incrementally using the Newton--Raphson method. At iteration $i$ of load step $n$, the displacement correction is obtained from
\begin{equation}
\mathbf{K}_t^{(i)}
\Delta \mathbf{d}^{(i)}
=
\mathbf{R}^{(i)},
\end{equation}
where $\mathbf{K}_t^{(i)}$ is the global tangent stiffness matrix assembled from the element tangent stiffness matrices. The iteration continues until the normalized residual satisfies
\begin{equation}
\frac{\|\mathbf{R}^{(i)}\|}{\|\mathbf{F}_{\mathrm{ext}}\|}
\le
\varepsilon_{\mathrm{tol}} .
\end{equation}

This formulation provides a unified nonlinear constitutive framework for polygonal CS-FEM. By selecting the appropriate tangent constitutive matrix, the proposed method can handle both elasto-plastic soil behavior and nonlinear stress-dependent rockfill behavior while preserving the advantages of strain smoothing and polygonal discretization.

\section{Implementation}
\label{sec:Implementation}

\subsection{Implementation of the elasto-plastic CS-FEM}

The proposed elasto-plastic CS-FEM is implemented in ABAQUS through the user element (UEL) interface. This implementation enables the polygonal CS-FEM formulation to be embedded into the incremental-iterative solution framework of ABAQUS/Standard. The UEL subroutine is responsible for evaluating and updating the element-level contributions to the residual force vector (RHS) and tangent stiffness matrix (AMATRX) using the nodal coordinates, element connectivity, displacement increments, material parameters, and state variables supplied by the ABAQUS solver~\citep{abaqus2010user}. Within this framework, the element tangent stiffness matrix and residual force vector are formulated as
\begin{equation}
\mathrm{AMATRX} = \tilde{\mathrm{K}}, \label{eq:amatrx}
\end{equation}
\begin{equation}
\mathrm{RHS} = \mathbf{F}_{\mathrm{ext}}^e - \mathbf{F}_{\mathrm{int}}^e, \label{eq:rhs}
\end{equation}
where $\mathbf{K}_t^e$ is the element tangent stiffness matrix, $\mathbf{F}_{\mathrm{ext}}^e$ is the external force vector, and $\mathbf{F}_{\mathrm{int}}^e$ is the internal force vector evaluated from the current stress state. The UEL routine iteratively updates both AMATRX and RHS according to Eqs.~(\ref{eq:amatrx}) and~(\ref{eq:rhs}) throughout the solution process.

The solution algorithm is summarized in Algorithm~\ref{alg1}. For each increment, the UEL subroutine computes the element-level stiffness matrices and residual forces by evaluating the Wachspress shape functions, the smoothed strain–displacement matrices, and the elasto-plastic constitutive matrices. These contributions are then assembled into the global tangent stiffness matrix and residual vector, which are used by ABAQUS to update the nodal displacements.

\begin{algorithm}[H]
\caption{Solving the nonlinear polygonal CS-FEM}
\begin{algorithmic}[1]
\Require Node and element information, material properties, nodal displacement $\mathbf{u}_t$
\Ensure Nodal displacement $\mathbf{u}_{t+1}$
\State \textbf{Iteration:} increment steps $k = 1$
\While{ABAQUS not converged}
    \State Solve nodal displacement $\mathbf{u}^{k}_{t+1}$
    \For{$i = 1$ to $n_c$}
        \State Solve shape functions $N_k(x)$ in Eq.~(\ref{eq:shape_function_values_2d})
        \State Solve strain–displacement matrix $\tilde{\mathbf{B}}_c$ in Eqs.~(\ref{eq:bk_2d})
        \State Solve elasto-plastic constitutive matrix $\mathbf{D}_{ep}$ in Eq.~(\ref{eq:CSFEM_Dep})
        \State Solve stiffness matrix of sub-cells $\tilde{\mathbf{K}}_c$ in Eq.~(\ref{eq:kmatrix})
    \EndFor
    \State Obtain element matrices $\tilde{\mathbf{K}}$ 
    \State Update AMATRX and RHS according to Eqs.~(\ref{eq:amatrx}) and (\ref{eq:rhs})
    \State $k = k + 1$
\EndWhile
\State Solve nodal field variables $\mathbf{u}_{t+1} = \mathbf{u}^{k}_{t+1}$
\end{algorithmic}
\label{alg1}
\end{algorithm}

\subsection{Post-processing of stress and strain}
\label{sec:stress_recovery}

After the element-level computations are completed via the UEL subroutine, the stress and strain tensors at the integration points are exported using the \texttt{UEXTERNALDB} interface in ABAQUS. This procedure provides the discrete material response at each Gauss point for all polygonal elements.

To obtain nodal stress and strain values suitable for visualization and further analysis, a stress recovery procedure is employed. In this approach, the nodal quantities are computed as weighted averages of the surrounding integration-point values:
\begin{equation}
\boldsymbol{\sigma}_i = \sum_{e \in \mathcal{E}_i} w_{i}^e \, \boldsymbol{\sigma}^e_\mathrm{gp}, 
\quad
\boldsymbol{\varepsilon}_i = \sum_{e \in \mathcal{E}_i} w_{i}^e \, \boldsymbol{\varepsilon}^e_\mathrm{gp},
\end{equation}
where $\mathcal{E}_i$ denotes the set of elements sharing node $i$, $\boldsymbol{\sigma}^e_\mathrm{gp}$ and $\boldsymbol{\varepsilon}^e_\mathrm{gp}$ are the stress and strain tensors at the integration points of element $e$, and $w_i^e$ are area- or volume-based weights satisfying the partition of unity condition:
\begin{equation}
\sum_{e \in \mathcal{E}_i} w_i^e = 1.
\end{equation}

This recovery ensures a continuous nodal stress and strain field from the discrete integration-point data. Such a continuous representation is essential for visualizing the stress and strain distributions over the polygonal mesh, for quantitatively evaluating peak and local stresses, and for performing post-processing tasks such as error estimation or verification against benchmark solutions.

\begin{figure}[H]
  \centering
  \includegraphics[width=1.0\textwidth]{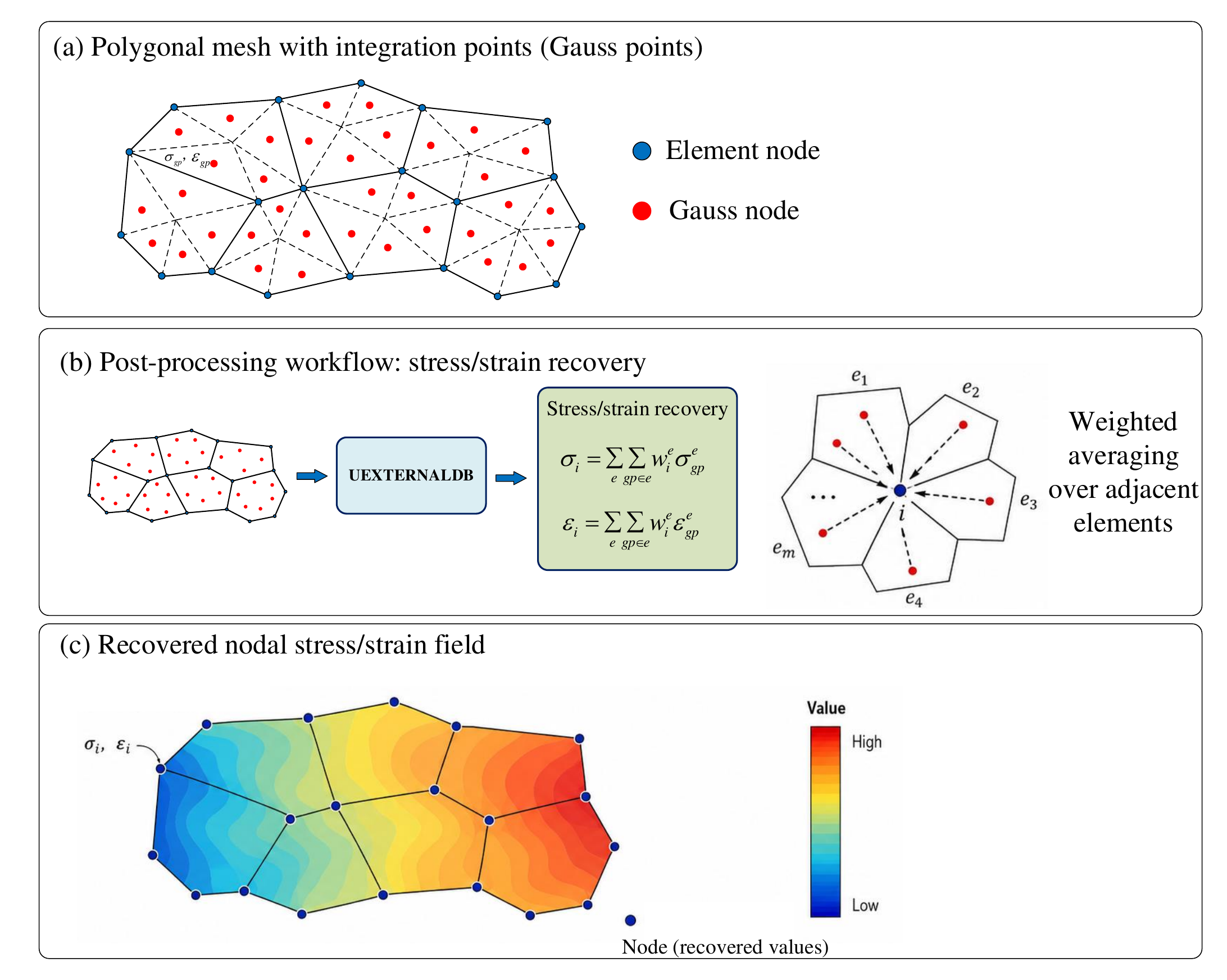}
  \caption{Schematic illustration of stress and strain recovery in polygonal CS-FEM. (a) Polygonal mesh with integration points (Gauss points) and nodes shared by neighboring elements. (b) Post-processing workflow using UEXTERNALDB and weighted averaging to recover nodal stress and strain. (c) Recovered nodal stress/strain field visualized over the polygonal mesh.}
  \label{fig:stress_recovery}
\end{figure}

A schematic workflow is illustrated in Fig.~\ref{fig:stress_recovery}, showing the mapping from integration-point data exported via \texttt{UEXTERNALDB} to smoothed nodal quantities.

\section{Numerical examples}
\label{sec:numerical_example}
\subsection{Uniaxial tension of an infinite plate with a circular hole}

\begin{figure}[H]
    \centering
    \includegraphics[width=1.0\linewidth]{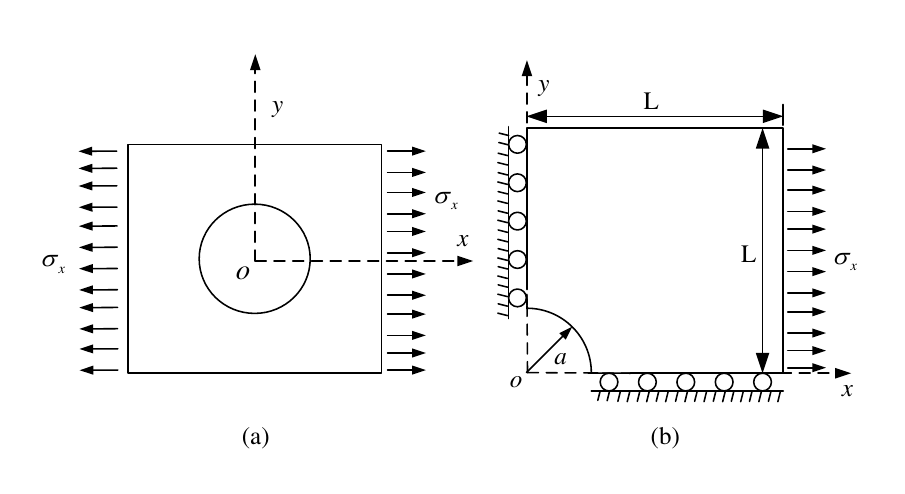}
    \caption{Geometry and boundary conditions for an infinite plate with a hole: (a) full infinite plate with a hole; (b) quarter model of the plate with a hole.}
    \label{fig:ex01_geo}
\end{figure}

An infinite plate with a circular hole of radius $a = 0.4\,\mathrm{m}$ is considered, as shown in Fig.~\ref{fig:ex01_geo}(a). The plate is subjected to a remote uniaxial tensile stress $\sigma_x = 1\,\mathrm{kPa}$ in the $x$ direction. Owing to geometric symmetry, only one quarter of the plate is modeled, as shown in Fig.~\ref{fig:ex01_geo}(b). The length of the quarter model is $L = 1.0\,\mathrm{m}$. The elastic modulus is $E = 1\times10^{5}\,\mathrm{Pa}$ and Poisson's ratio is $\nu = 0.25$. Normal displacement constraints are imposed on the bottom boundary and the left boundary of the quarter plate.

The analytical solution in polar coordinates $(r,\theta)$ is given by \citep{Song2018SBFEM}
\begin{equation}
u_x=\frac{Pa}{8G}\left[\frac{r}{a}(1+\kappa)\cos\theta+\frac{2a}{r}\left((1+\kappa)\cos\theta+\cos3\theta\right)-\frac{2a^3}{r^3}\cos3\theta\right],
\end{equation}
\begin{equation}
u_y=\frac{Pa}{8G}\left[\frac{r}{a}(\kappa-3)\sin\theta+\frac{2a}{r}\left((1-\kappa)\sin\theta+\sin3\theta\right)-\frac{2a^3}{r^3}\sin3\theta\right],
\end{equation}
where the shear modulus $G$ and the Kolosov constant $\kappa$ are defined as
\begin{equation}
G=\frac{E}{2(1+\nu)},
\end{equation}
\begin{equation}
\kappa=\frac{3-\nu}{1+\nu}.
\end{equation}

\begin{figure}[H]
    \centering
    \includegraphics[width=1.0\linewidth]{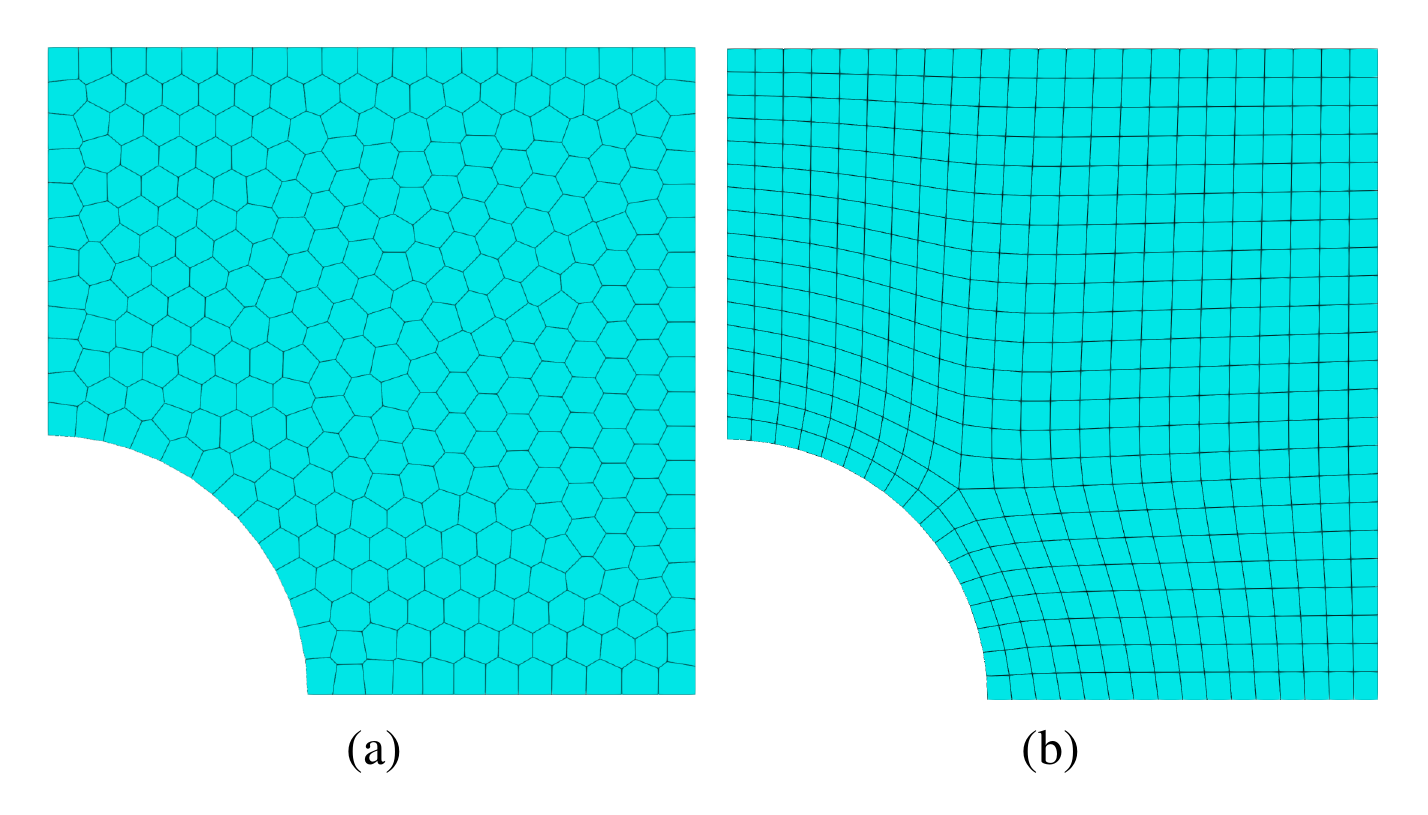}
    \caption{Mesh models of the infinite plate with a hole (element size: 0.05 m): (a) polygonal mesh; (b) quadrilateral mesh.}
    \label{fig:ex01_mesh}
\end{figure}

To verify the accuracy and convergence of the proposed method, four different element sizes, namely $0.1\,\mathrm{m}$, $0.05\,\mathrm{m}$, $0.025\,\mathrm{m}$, and $0.0125\,\mathrm{m}$, are considered. For the CS-FEM model, polygonal meshes are employed, as shown in Fig.~\ref{fig:ex01_mesh}(a), and plane stress elements are used. For the conventional FEM model, quadrilateral meshes are adopted, as shown in Fig.~\ref{fig:ex01_mesh}(b). ABAQUS is used for comparison, with CPS4 elements adopted in the FEM simulations.

\begin{figure}[H]
    \centering
    \includegraphics[width=1.0\linewidth]{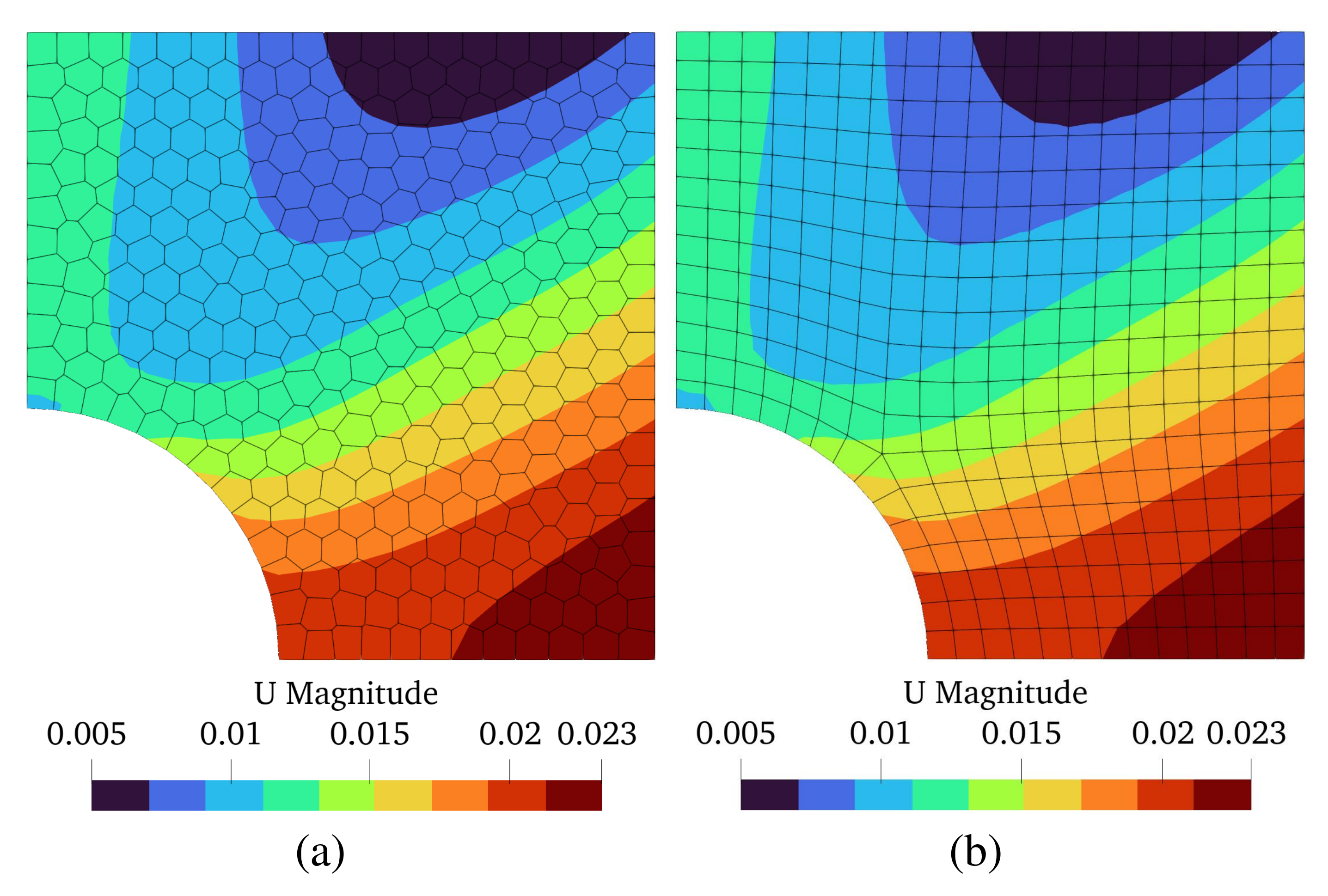}
    \caption{Displacement contours of the infinite plate with a hole (element size: 0.05 m): (a) CS-FEM; (b) FEM.}
    \label{fig:ex01_contour}
\end{figure}

\begin{figure}[H]
    \centering
    \includegraphics[width=0.7\linewidth]{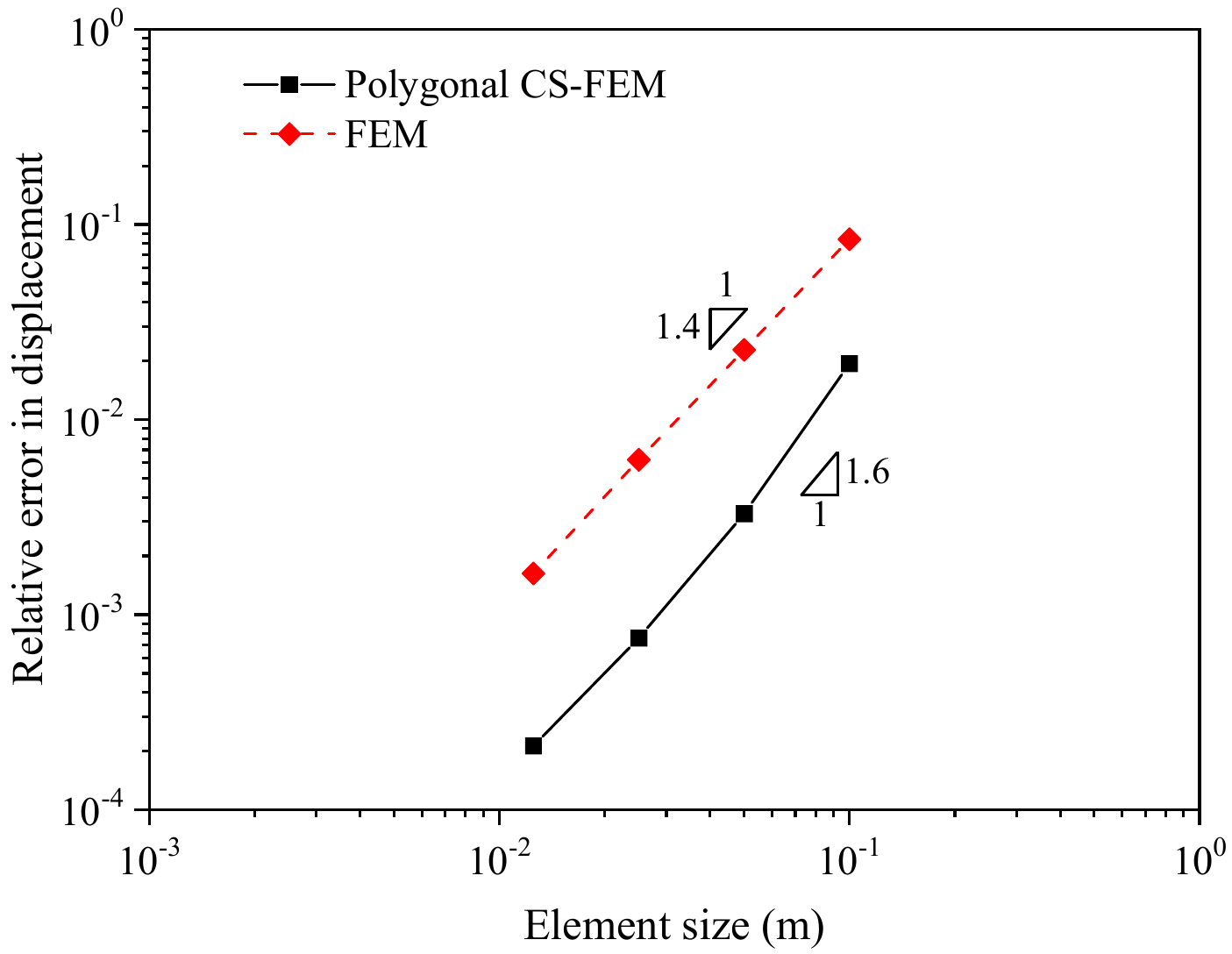}
    \caption{Convergence of the relative displacement error for the infinite plate with a circular hole.}
    \label{fig:ex01_error}
\end{figure}

Fig.~\ref{fig:ex01_contour} shows the displacement contours of the perforated plate for the element size of $0.05\,\mathrm{m}$. It can be seen that the displacement distributions obtained by polygonal CS-FEM and FEM are generally consistent. In addition, the displacement field predicted by CS-FEM is smoother than that obtained by FEM. Fig.~\ref{fig:ex01_error} presents the convergence of the relative displacement error for the infinite plate with a circular hole. The relative errors of both methods decrease with mesh refinement. The estimated convergence slopes are approximately 1.6 for the polygonal CS-FEM and 1.4 for the FEM, showing that the polygonal CS-FEM achieves a slightly higher convergence rate. In addition, for the same element size, the displacement error of the polygonal CS-FEM is consistently lower than that of the FEM. This confirms the improved accuracy of the proposed method.

\subsection{Bearing capacity problem}

\begin{figure}[H]
  \centering
  \includegraphics[width=0.9\textwidth]{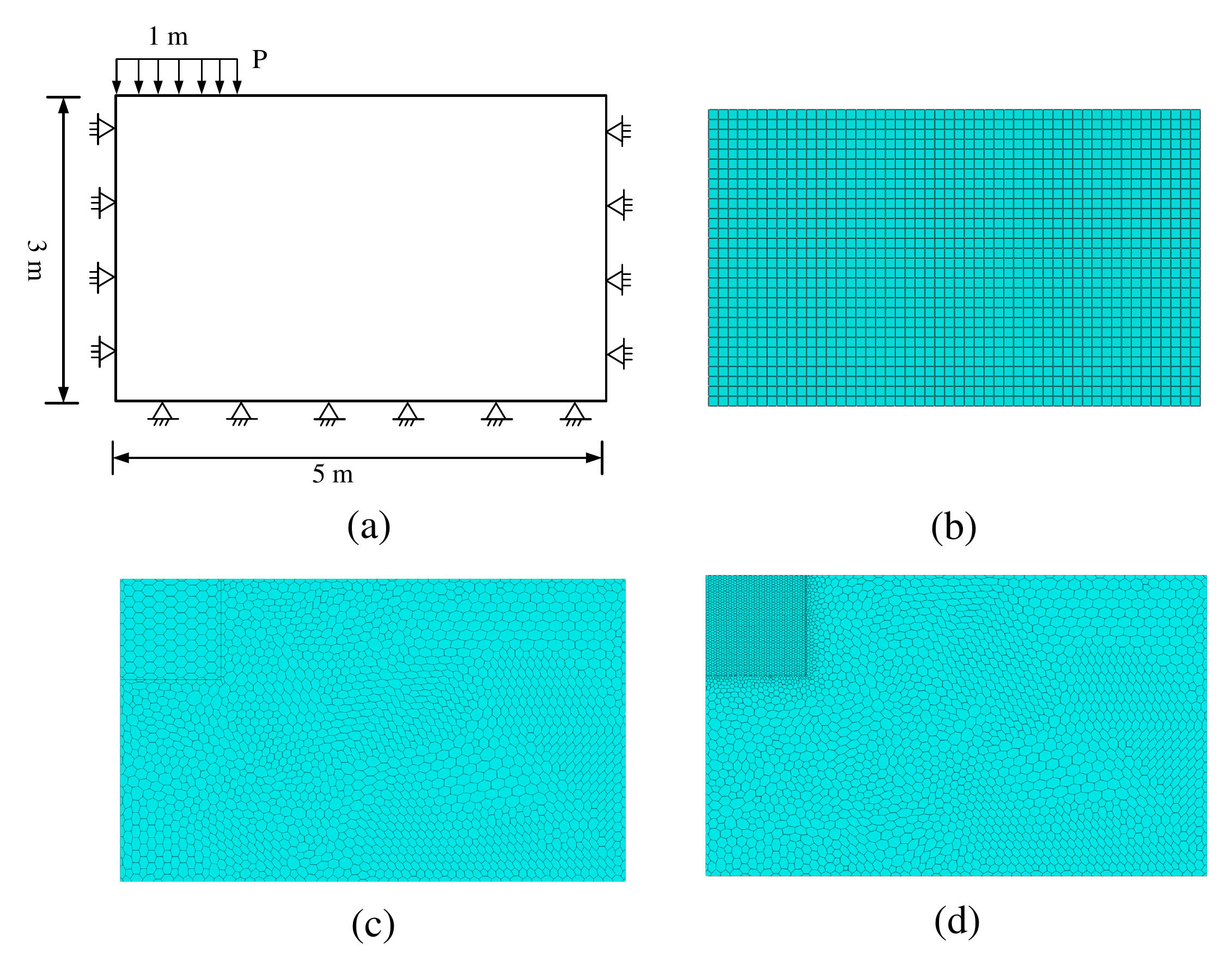}
  \caption{{Schematic diagram of the bearing capacity problem: (a) geometry and boundary conditions; (b) quadrilateral mesh; (c) polygonal mesh; (d) locally refined polygonal mesh.}}
  \label{fig:ex02_geo_mesh}
\end{figure}

In this example, a flexible strip footing resting on the surface of a half-space soil domain is considered, as shown in Fig.~\ref{fig:ex02_geo_mesh}(a). Due to symmetry, only half of the footing is modeled. Roller boundary conditions are imposed on the two side boundaries, while the bottom boundary is fully fixed. A uniformly distributed vertical load is applied over a footing width of 1~m. The computational domain has a width of 5~m and a height of 3~m. Three types of meshes are used in this example, including a quadrilateral mesh for FEM, a polygonal mesh for polygonal CS-FEM, and a locally refined polygonal mesh for polygonal CS-FEM, as shown in Figs.~\ref{fig:ex02_geo_mesh}(b)--(d). The corresponding soil parameters are listed in Tab.~\ref{tab:ex02:soil_parameters}. For a soil with zero unit weight ($\gamma = 0$), the ultimate bearing capacity of strip foundations can be expressed as \citep{sherif2013matlab}:

\begin{equation}
q_u=\begin{cases}
c\cdot N_c & \phi \neq 0, \\
(\pi+2)\cdot c & \phi = 0,
\end{cases}
\label{eq:q}
\end{equation}
where $q_u$ is the ultimate bearing capacity of the strip footing, 
$c$ is the soil cohesion, $\phi$ is the internal friction angle, and 
$N_c$ is the bearing capacity factor. For $\phi=0$, the expression reduces 
to the classical solution for purely cohesive soil. The bearing capacity 
factor $N_c$ is given by
\begin{equation}
N_c=(N_q-1)\cot\phi,
\end{equation}
where $N_q$ is another bearing capacity factor defined as
\begin{equation}
N_q=e^{\pi\tan\phi}\left(\frac{1+\sin\phi}{1-\sin\phi}\right).
\end{equation}

\begin{table}[H]
  \centering
  \caption{Soil parameters of a flexible strip footing.}
  \label{tab:ex02:soil_parameters}
  \begin{tabular}{ccccc}
    \toprule
    $E$ (kPa) & $\nu$ & $\phi$ & $c$ (kPa) & $\psi$ \\ 
    \midrule
    10000 & 0.3 & $5^{\circ}$ & 1 & $5^{\circ}$ \\ 
    \bottomrule
  \end{tabular}
  
  \vspace{2mm}
  {\footnotesize Note: $\psi$ denotes the dilation angle.}
\end{table}

\begin{figure}[H]
  \centering
  \includegraphics[width=0.9\textwidth]{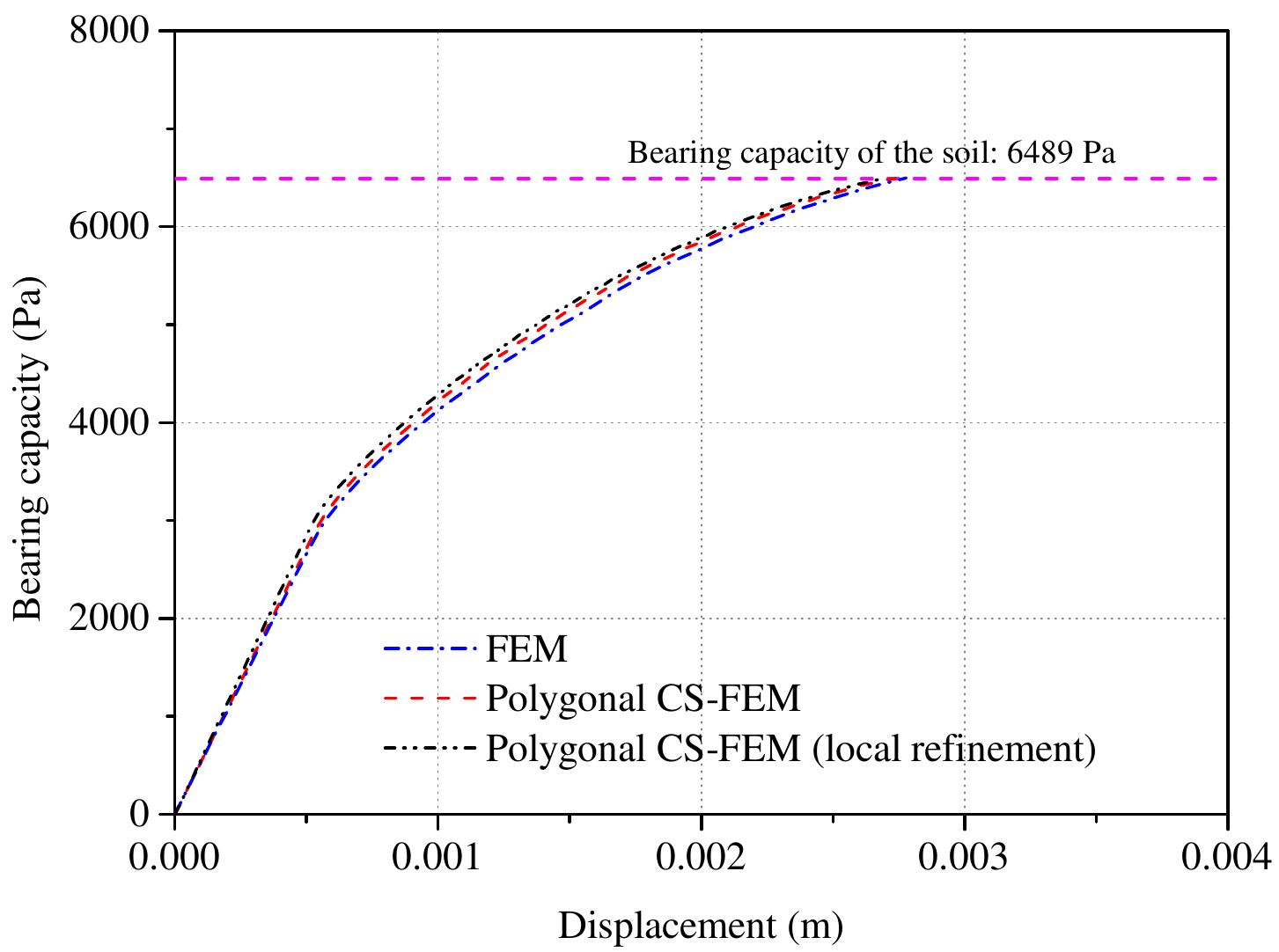}
  \caption{\textcolor{black}{Bearing capacity--displacement diagram of a flexible footing using polygonal CS-FEM and FEM.}}
  \label{fig:ex02_disp_stress}
\end{figure}

Fig.~\ref{fig:ex02_disp_stress} presents the bearing capacity--displacement curves obtained by FEM, polygonal CS-FEM, and polygonal CS-FEM with local refinement. The three numerical results show very similar nonlinear response trends and gradually approach the analytical bearing capacity. According to Eq.~\eqref{eq:q}, the analytical bearing capacity of the soil is 6489~Pa. As shown in Tab.~\ref{tab:ex02_relative_errors}, the bearing capacities predicted by polygonal CS-FEM, polygonal CS-FEM with local refinement, and FEM are 6495.93~Pa, 6493.41~Pa, and 6499.53~Pa, respectively. The corresponding relative errors are $1.07\times10^{-3}$, $6.80\times10^{-4}$, and $1.62\times10^{-3}$, respectively. These results indicate that the polygonal CS-FEM provides more accurate predictions than the conventional FEM, and the locally refined polygonal mesh further improves the accuracy.

\begin{figure}[H]
  \centering
  \includegraphics[width=1.0\textwidth]{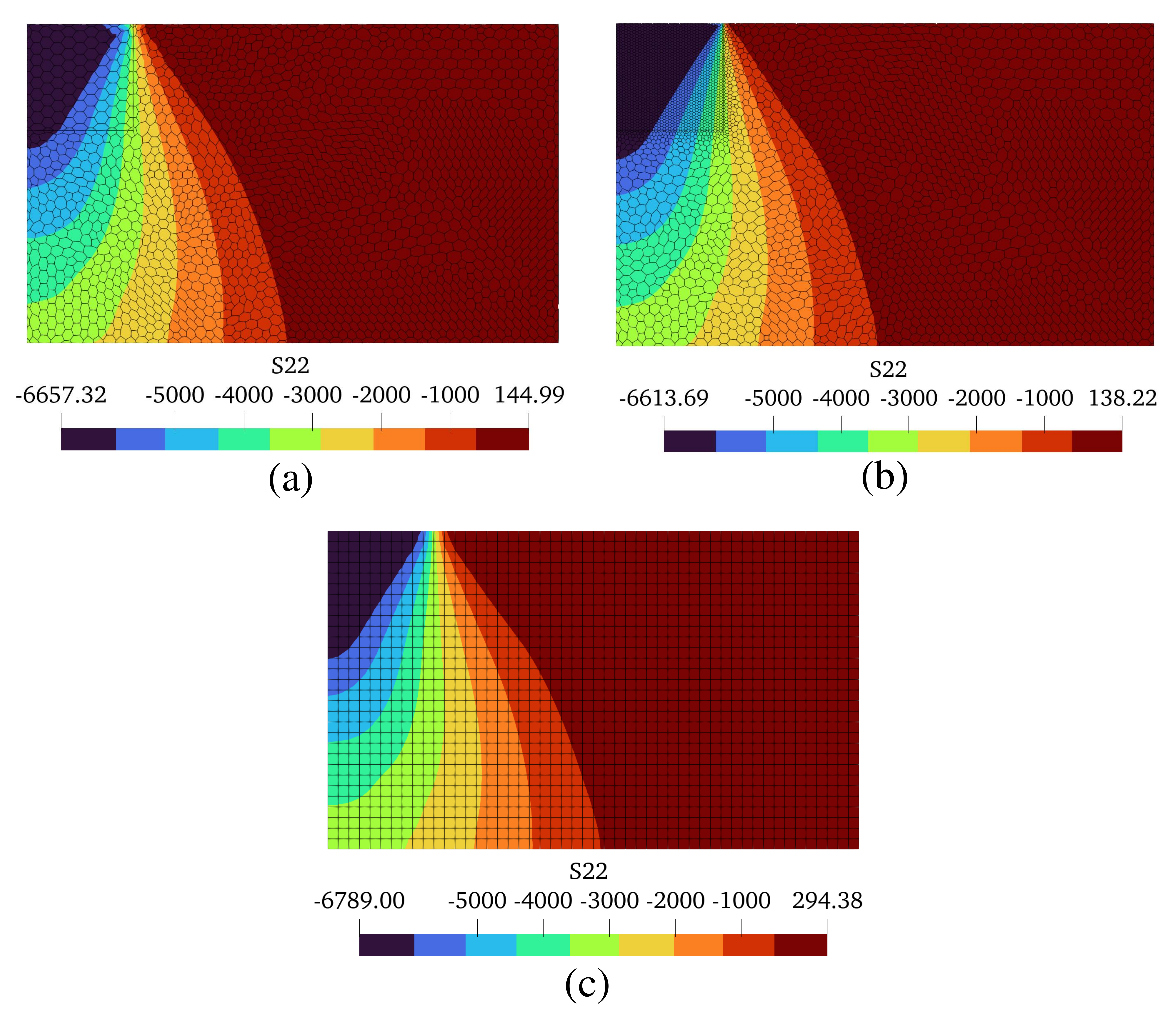}
  \caption{\textcolor{black}{Comparison of stress contours between polygonal CS-FEM and FEM at the ultimate bearing capacity: (a) polygonal CS-FEM; (b) polygonal CS-FEM with local refinement; (c) FEM.}}
  \label{fig:ex02_contour}
\end{figure}

Fig.~\ref{fig:ex02_contour} compares the stress contours obtained by polygonal CS-FEM, polygonal CS-FEM with local refinement, and FEM at the ultimate bearing capacity. The stress concentration zone beneath the footing can be clearly captured by all three methods, and the overall stress distributions are generally consistent with each other. In the locally refined polygonal CS-FEM result, the stress variation near the loaded region is represented with more local details due to the refined discretization. Together with the bearing-capacity errors listed in Tab.~\ref{tab:ex02_relative_errors}, these results indicate that the polygonal CS-FEM can accurately reproduce the stress response of the footing problem while providing improved bearing-capacity prediction compared with the conventional FEM.

\begin{table}[H]
  \centering
  \caption{{Relative errors in bearing capacity between polygonal CS-FEM and FEM.}}
  \label{tab:ex02_relative_errors}
  \resizebox{\textwidth}{!}{
  \begin{tabular}{ccccccc}
    \toprule
    \multirow{3}{*}{\shortstack{Analytical\\solution (Pa)}}
    & \multicolumn{3}{c}{Bearing capacity (Pa)}
    & \multicolumn{3}{c}{Relative error} \\
    \cmidrule(lr){2-4} \cmidrule(lr){5-7}
    &
    \shortstack{Polygonal\\CS-FEM}
    & \shortstack{Polygonal CS-FEM\\(local refinement)}
    & FEM
    & \shortstack{Polygonal\\CS-FEM}
    & \shortstack{Polygonal CS-FEM\\(local refinement)}
    & FEM \\
    \midrule
    6489.00
    & 6495.93
    & 6493.41
    & 6499.53
    & $1.07\times10^{-3}$
    & $6.80\times10^{-4}$
    & $1.62\times10^{-3}$ \\
    \bottomrule
  \end{tabular}
  }
\end{table}

\subsection{Staged construction analysis of a core rockfill dam}

\begin{figure}[H]
    \centering
    \includegraphics[width=1.0\linewidth]{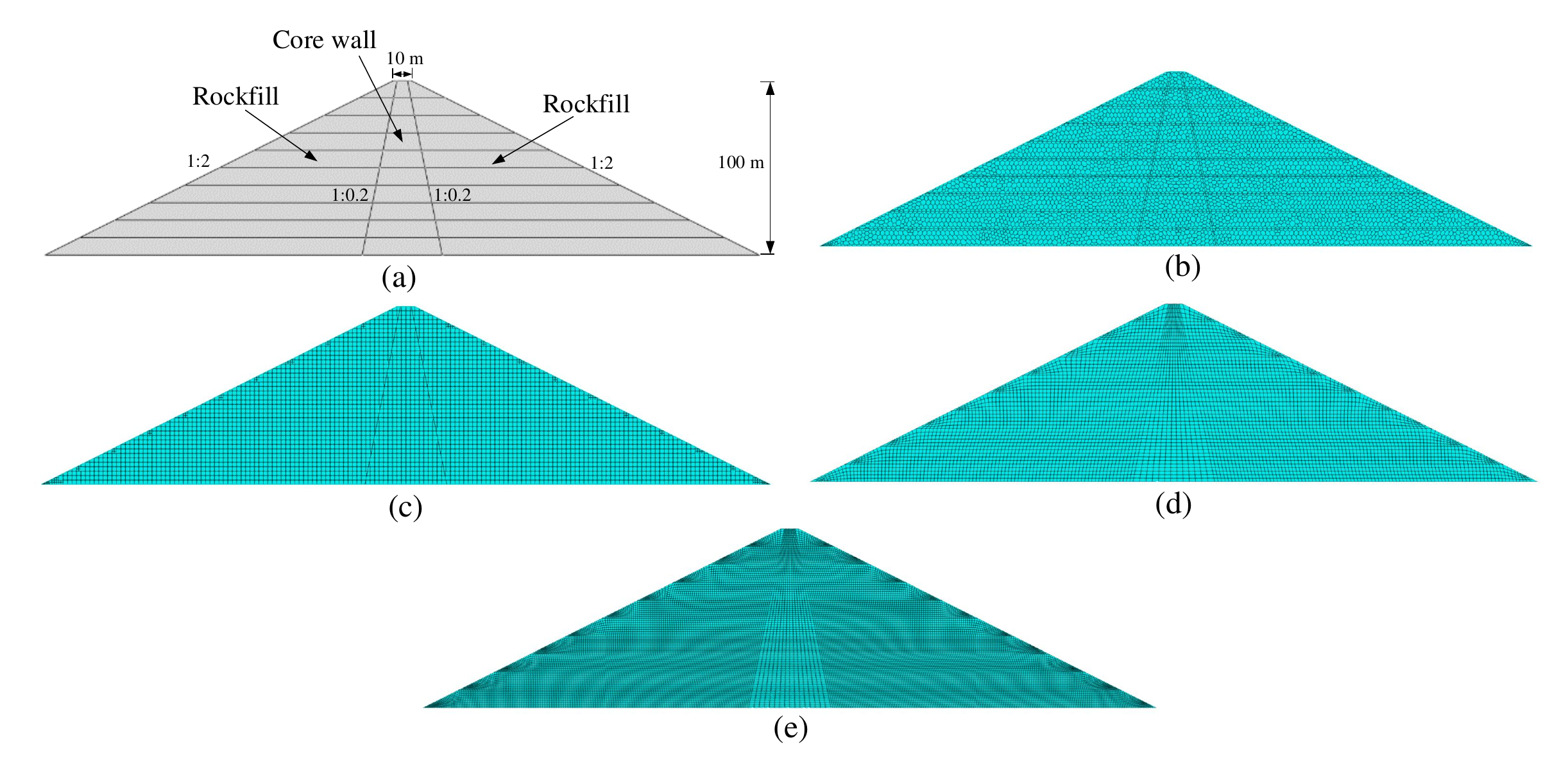}
    \caption{Geometry and mesh models of the rockfill dam: (a) geometry model; (b) polygonal mesh; (c) hybrid quadtree mesh; (d) FEM mesh; and (e) refined FEM mesh used for the reference solution.}
    \label{fig:ex03_geo_mesh}
\end{figure}

An idealized core rockfill dam is analyzed to evaluate the performance of the proposed method in staged construction problems. As shown in Fig.~\ref{fig:ex03_geo_mesh}(a), the dam is 100 m high with a crest width of 10 m, and both the upstream and downstream slopes are 1:2. The central core wall has a crest width of 6 m and side slopes of 1:0.2. The dam construction is simulated using 10 filling layers, each with a thickness of 10 m, so that the stress evolution during construction can be captured.

The computational models are shown in Fig.~\ref{fig:ex03_geo_mesh}. Two polygonal CS-FEM discretizations are considered, including an irregular polygonal mesh and a hybrid quadtree mesh, both with a characteristic element size of 2 m. A conventional FEM model with quadrilateral elements and the same characteristic element size is also established for comparison. In addition, a refined FEM mesh with an element size of 1 m is used to provide the reference solution. The dam materials are modeled using the Duncan--Chang (E--B) model, and the corresponding parameters are listed in Tab.~\ref{tab:DC_mat}.

\begin{figure}[H]
    \centering
    \includegraphics[width=1.0\linewidth]{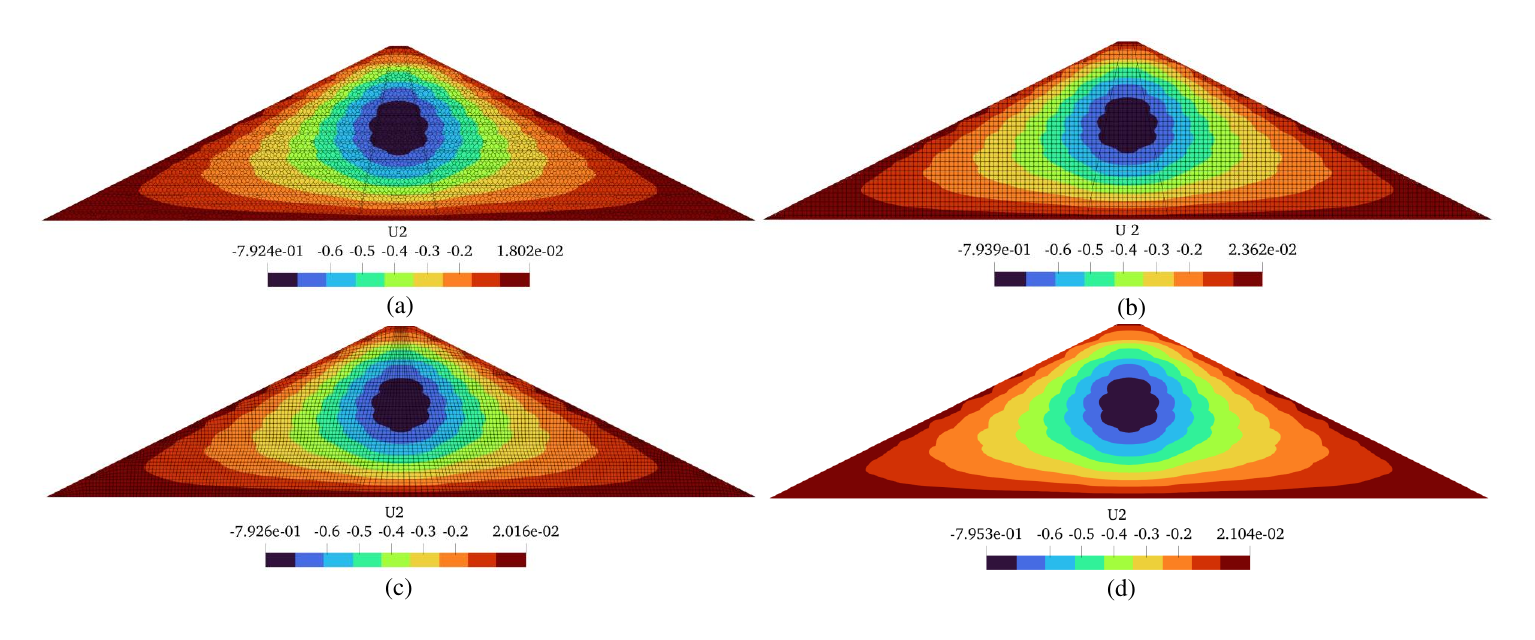}
    \caption{Comparison of dam settlement distributions obtained by different numerical models: (a) CS-FEM with polygonal elements; (b) CS-FEM with quadtree elements; (c) FEM result; and (d) reference solution.}
    \label{fig:ex03_disp1}
\end{figure}

\begin{figure}[H]
    \centering
    \includegraphics[width=1.0\linewidth]{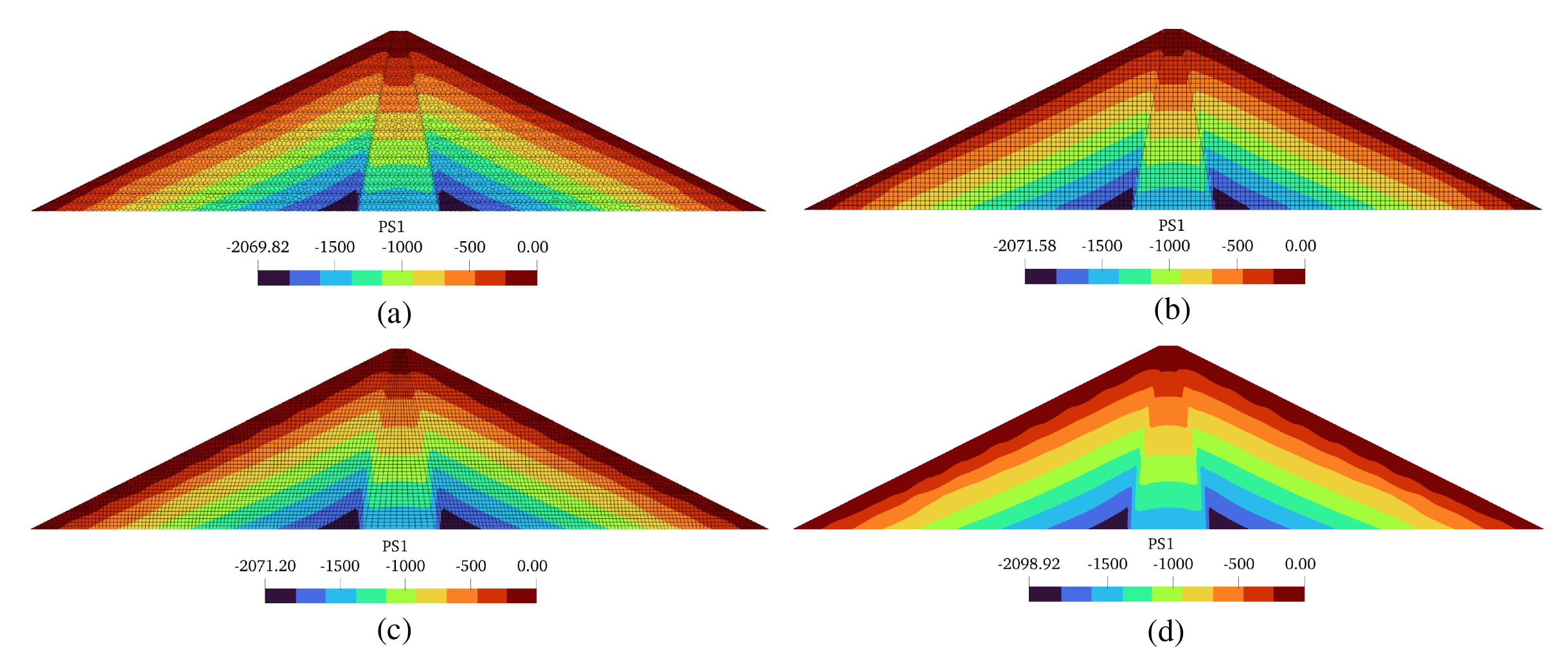}
    \caption{Comparison of major principal stress distributions in the dam obtained by different numerical models: (a) CS-FEM with polygonal elements; (b) CS-FEM with quadtree elements; (c) FEM result; and (d) reference solution.}
    \label{fig:ex03_ps}
\end{figure}

Figs.~\ref{fig:ex03_disp1} and~\ref{fig:ex03_ps} compare the settlement and major principal stress distributions obtained from different numerical models. All methods produce similar deformation and stress patterns. The maximum settlement occurs near the dam axis and core wall, while the major principal stress gradually increases with depth and shows a clear layered distribution caused by staged filling.

\begin{table}[H]
\centering
\caption{Duncan--Chang (E--B) model parameters for the core wall and rockfill materials}
\resizebox{\textwidth}{!}{
\begin{tabular}{lcccccccccc}
\hline
Material & $K$ & $n$ & $R_f$ & $c$ (kPa) & $\varphi$ ($^\circ$) & $\varphi_0$ ($^\circ$) & $K_{ur}$ & $K_b$ & $m$ & $\rho$ (g/cm$^3$) \\
\hline
Core wall      & 500  & 0.35 & 0.8 & 50 & 30 & 0 & 800  & 470 & 0.15 & 2.0 \\
Rockfill       & 1100 & 0.30 & 0.8 & 10 & 40 & 0 & 1800 & 600 & 0.10 & 2.2 \\
\hline
\end{tabular}
}
\label{tab:DC_mat}
\end{table}

The quantitative comparison is given in Tab.~\ref{tab:dam_settlement_stress_comparison}. The settlement errors of the polygonal CS-FEM, quadtree CS-FEM, and FEM are $3.65\times10^{-3}$, $1.76\times10^{-3}$, and $3.40\times10^{-3}$, respectively. The corresponding stress errors are $1.39\times10^{-2}$, $1.30\times10^{-2}$, and $1.32\times10^{-2}$. These results indicate that the proposed CS-FEM gives results close to the reference solution, and the quadtree discretization provides slightly better settlement accuracy.

\begin{figure}[H]
    \centering
    \includegraphics[width=0.8\linewidth]{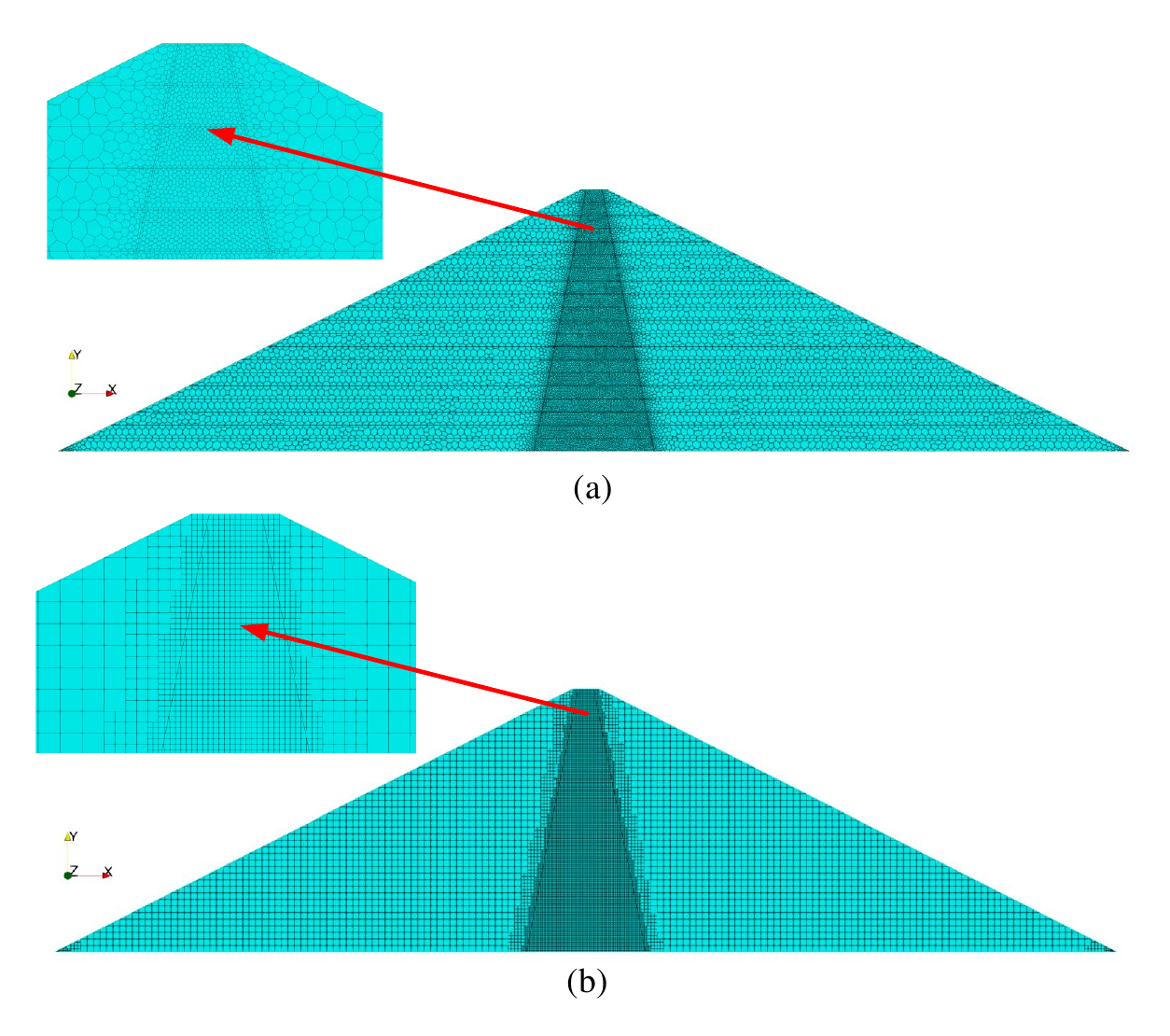}
    \caption{Comparison of mesh refinement strategies for the rockfill dam with core wall refinement: (a) polygonal mesh enabling flexible local refinement in the core wall; (b) hybrid quadtree mesh with hierarchical refinement in the core wall.}
    \label{fig:ex03_mesh_jm}
\end{figure}

\begin{figure}[H]
    \centering
    \includegraphics[width=1.0\linewidth]{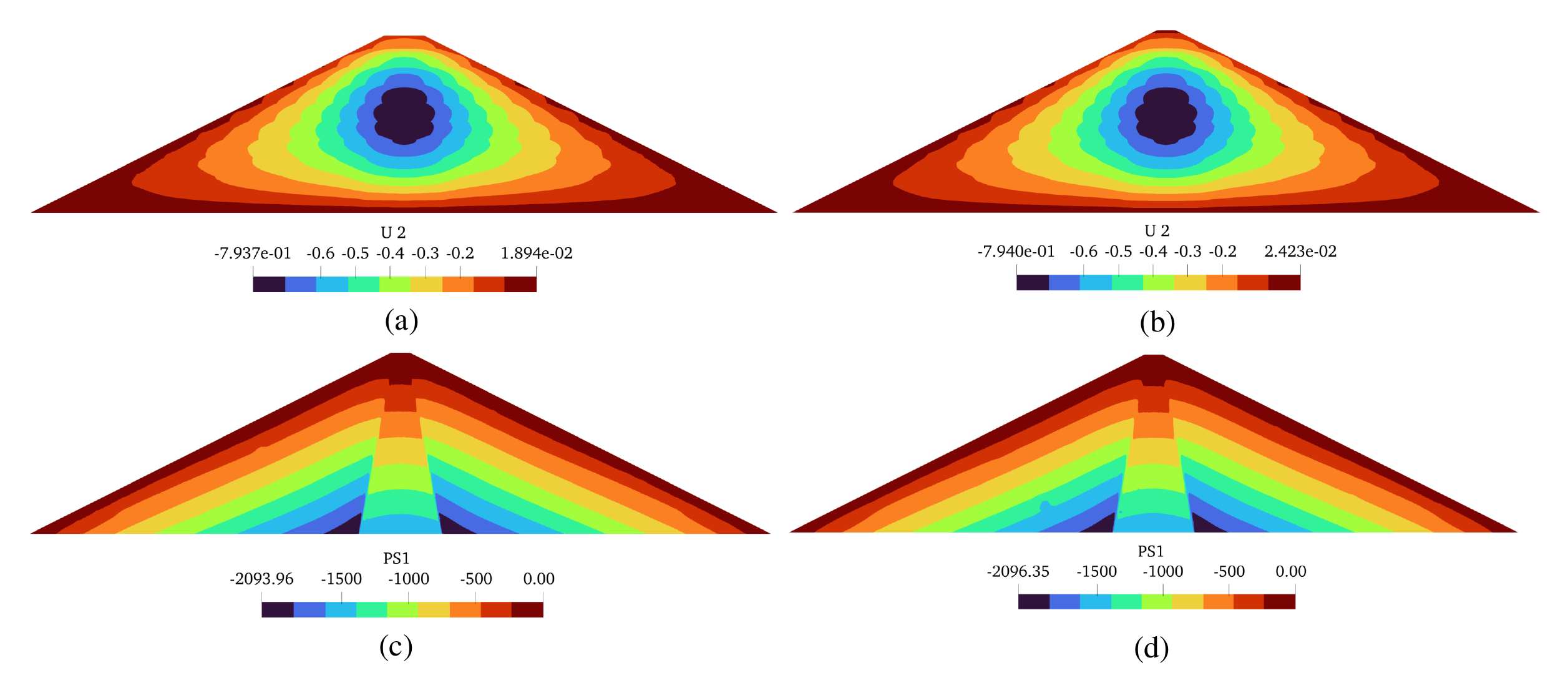}
    \caption{Comparison of settlement and major principal stress distributions after local refinement of the core wall: (a) displacement using polygonal elements; (b) displacement using quadtree elements; (c) major principal stress using polygonal elements; (d) major principal stress using quadtree elements.}
    \label{fig:ex03_result_jm}
\end{figure}

To further examine the local refinement capability of the proposed method, the core wall region is refined using both polygonal and quadtree meshes, as shown in Fig.~\ref{fig:ex03_mesh_jm}. Fig.~\ref{fig:ex03_result_jm} shows the settlement and major principal stress distributions after local refinement of the core wall. Compared with the non-refined results, the refined models provide a more detailed representation of the local response near the core wall while preserving the overall deformation and stress patterns of the dam.

The errors after local refinement are summarized in Tab.~\ref{tab:core_refined_settlement_stress_comparison}. The settlement errors decrease to $2.01\times10^{-3}$ for the polygonal mesh and $1.63\times10^{-3}$ for the quadtree mesh. More importantly, the stress errors are reduced to $2.36\times10^{-3}$ and $1.22\times10^{-3}$, respectively. This confirms that local refinement of the core wall significantly improves the stress prediction accuracy, especially for the quadtree CS-FEM model. 

\begin{table}[H]
\centering
\caption{Comparison of dam settlement and major principal stress obtained by different numerical models.}
\label{tab:dam_settlement_stress_comparison}
\resizebox{\textwidth}{!}{
\begin{tabular}{ccccc}
\toprule
Model 
& Settlement (m) 
& Settlement error 
& Major principal stress (kPa) 
& Stress error \\
\midrule
CS-FEM with polygonal elements 
& $7.924\times10^{-1}$ 
& $3.65\times10^{-3}$ 
& $-2.06982\times10^{3}$ 
& $1.39\times10^{-2}$ \\

CS-FEM with quadtree elements  
& $7.939\times10^{-1}$ 
& $1.76\times10^{-3}$ 
& $-2.07158\times10^{3}$ 
& $1.30\times10^{-2}$ \\

FEM                           
& $7.926\times10^{-1}$ 
& $3.40\times10^{-3}$ 
& $-2.07120\times10^{3}$ 
& $1.32\times10^{-2}$ \\

Reference solution             
& $7.953\times10^{-1}$ 
& -- 
& $-2.09892\times10^{3}$ 
& -- \\
\bottomrule
\end{tabular}
}
\end{table}

\begin{table}[H]
\centering
\caption{Comparison of dam-body settlement and major principal stress after local refinement of the core wall.}
\label{tab:core_refined_settlement_stress_comparison}
\resizebox{\textwidth}{!}{
\begin{tabular}{ccccc}
\toprule
Model 
& Settlement (m) 
& Settlement error 
& Major principal stress (kPa) 
& Stress error \\
\midrule
CS-FEM with polygonal elements 
& $7.937\times10^{-1}$ 
& $2.01\times10^{-3}$ 
& $-2.09396\times10^{3}$ 
& $2.36\times10^{-3}$ \\

CS-FEM with quadtree elements  
& $7.940\times10^{-1}$ 
& $1.63\times10^{-3}$ 
& $-2.09635\times10^{3}$ 
& $1.22\times10^{-3}$ \\

Reference solution             
& $7.953\times10^{-1}$ 
& -- 
& $-2.09892\times10^{3}$ 
& -- \\
\bottomrule
\end{tabular}
}
\end{table}

\subsection{Tunnel excavation}
\begin{figure}[H]
    \centering
    \includegraphics[width=1.0\linewidth]{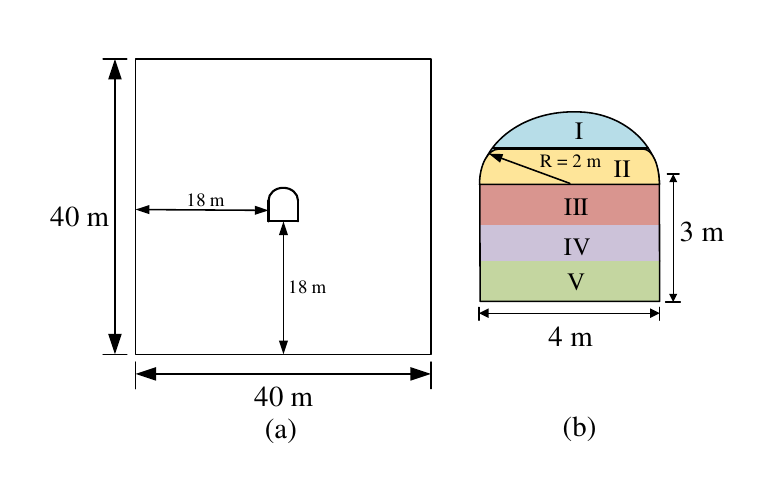}
\caption{Geometry and staged excavation sequence of the tunnel problem; (a) computational domain and tunnel location; (b) excavation sequence from stage I to stage V.}
\label{fig:ex04_geo}
\end{figure}

A tunnel excavation problem is analyzed in this example to further evaluate the applicability of the proposed method to staged excavation analysis. The computational domain and tunnel location are shown in Fig.~\ref{fig:ex04_geo}(a). The tunnel has a gate-shaped cross section, consisting of two vertical sidewalls and a semicircular crown arch. The tunnel width is 4~m, the sidewall height is 3~m, and the crown arch radius is 2~m. The excavation process is simulated using five excavation stages, as shown in Fig.~\ref{fig:ex04_geo}(b), with each stage corresponding to an excavation thickness of 1~m. The surrounding rock is assumed to follow the Mohr--Coulomb failure criterion. The elastic modulus, Poisson's ratio, cohesion, and internal friction angle are $E = 8~\mathrm{MPa}$, $\nu = 0.27$, $c = 110~\mathrm{kPa}$, and $\varphi = 39^\circ$, respectively, as summarized in Tab.~\ref{tab:ex04_rock_parameters}. 

\begin{figure}[H]
    \centering
    \includegraphics[width=1.0\linewidth]{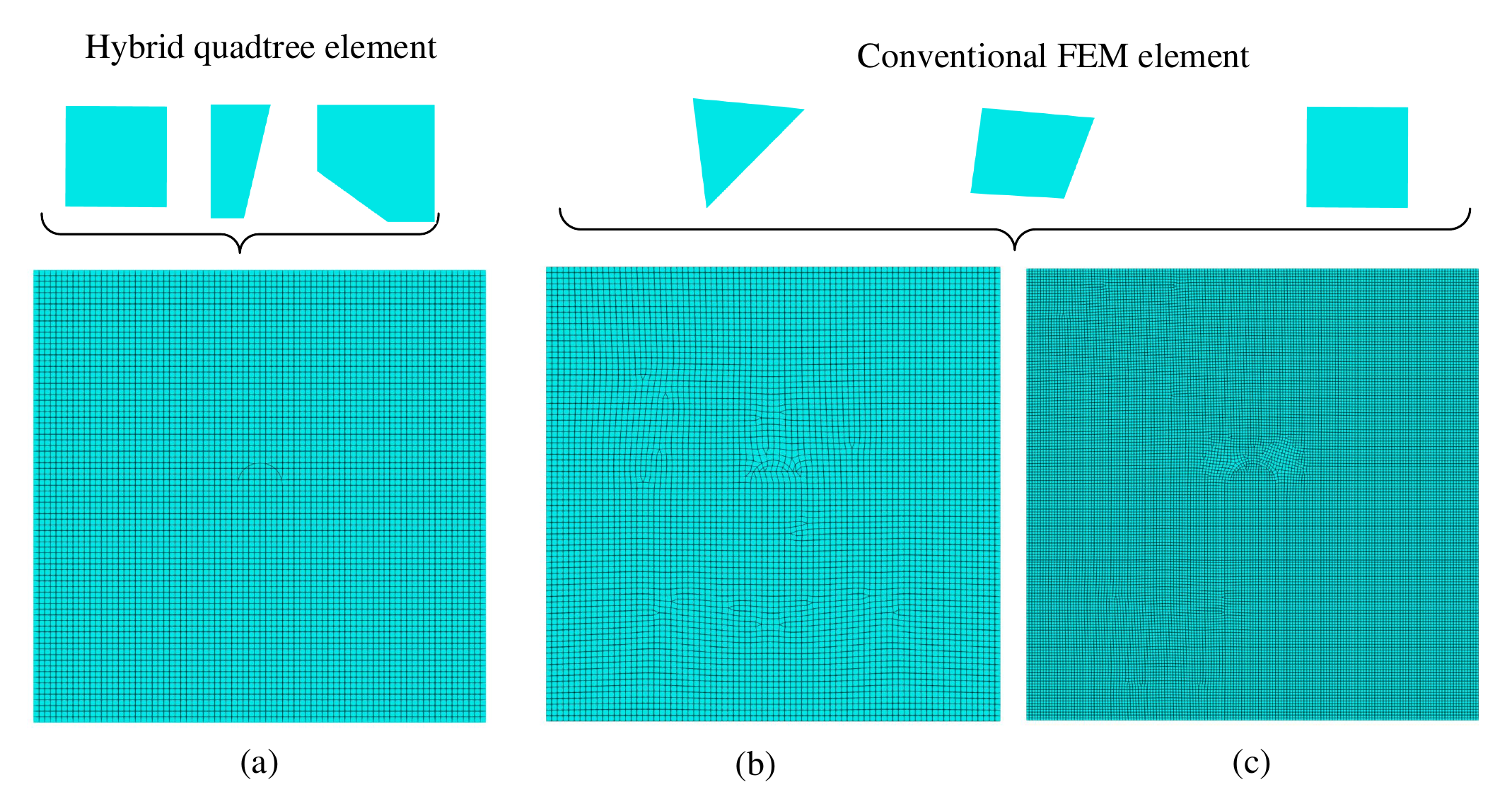}
    \caption{Mesh discretizations used for the tunnel excavation problem; (a) hybrid quadtree mesh for the proposed CS-FEM; (b) conventional FEM mesh; (c) refined FEM mesh for the reference solution.}
    \label{fig:ex04_mesh}
\end{figure}

Three types of mesh discretizations are used for comparison, as shown in Fig.~\ref{fig:ex04_mesh}. The hybrid quadtree mesh used in the proposed polygonal CS-FEM is shown in Fig.~\ref{fig:ex04_mesh}(a), with a characteristic mesh size of 0.5 m. This mesh consists of regular quadrilateral elements and other polygonal elements. Fig.~\ref{fig:ex04_mesh}(b) presents the conventional FEM mesh with the same mesh size of 0.5 m for comparison, while Fig.~\ref{fig:ex04_mesh}(c) shows a refined FEM mesh with a mesh size of 0.25 m, which is used to provide the reference solution.

\begin{figure}[H]
    \centering
    \includegraphics[width=1.0\linewidth]{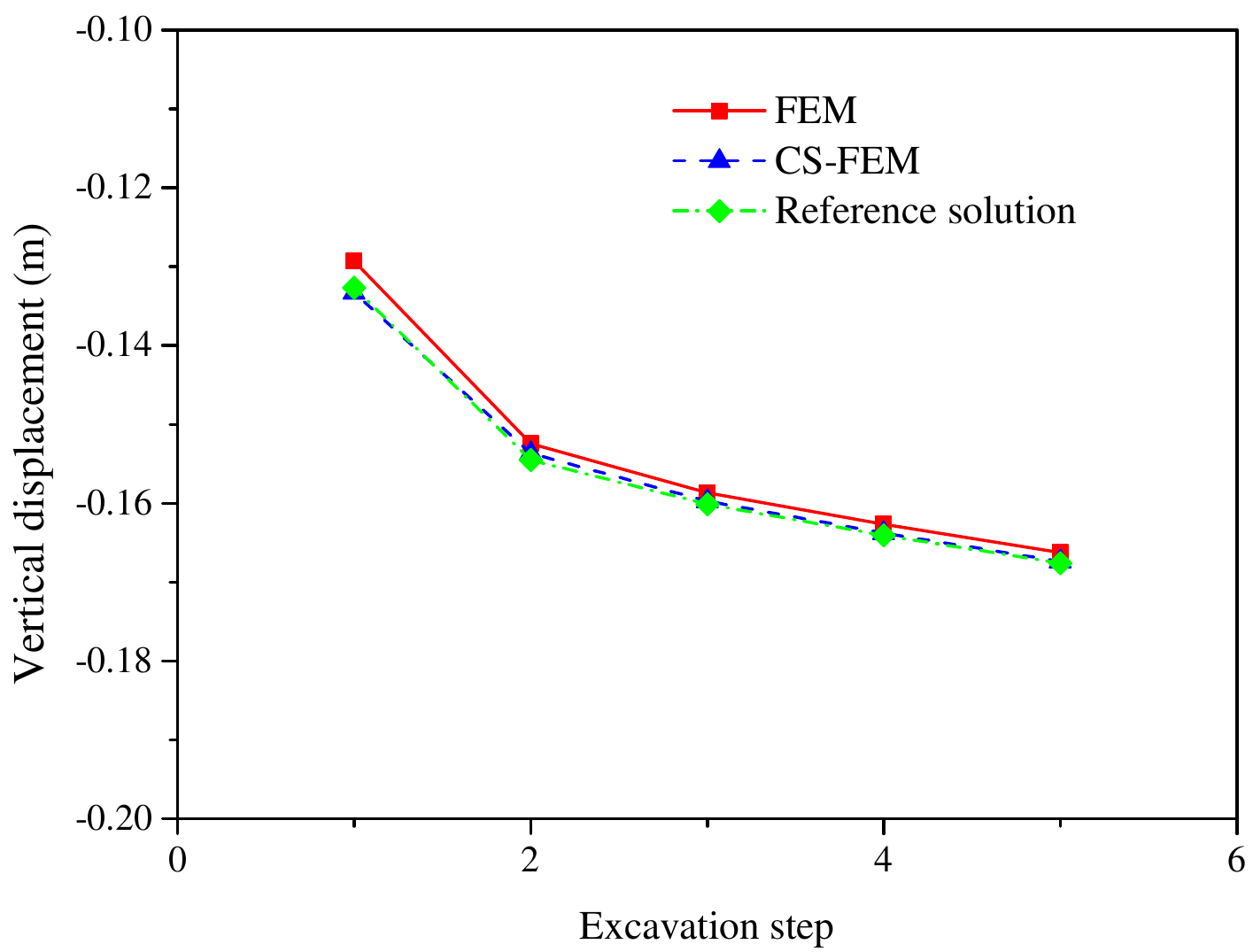}
    \caption{Comparison of vertical displacement at the tunnel crown obtained by different numerical methods under different excavation steps.}
    \label{fig:ex04_staged}
\end{figure}

\begin{figure}[H]
    \centering
    \includegraphics[width=1.0\linewidth]{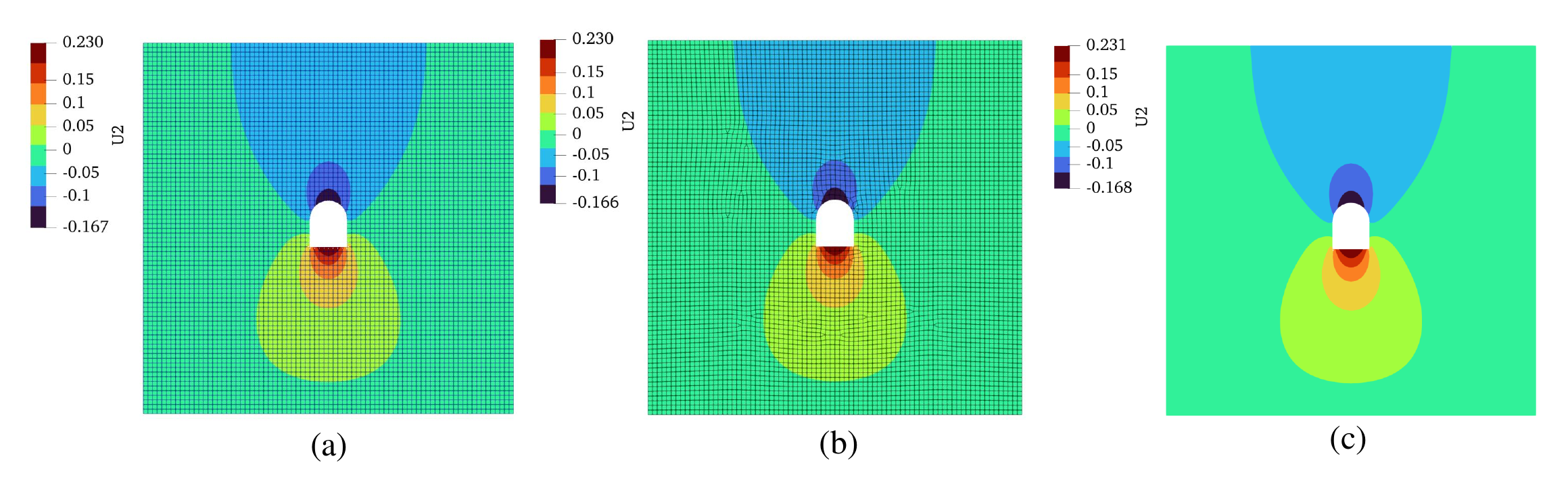}
\caption{Comparison of vertical displacement fields for the gate-shaped tunnel excavation problem: (a) polygonal CS-FEM result; (b) conventional FEM result;  (c) reference solution.}
    \label{fig:ex04_contour}
\end{figure}

As shown in Fig.~\ref{fig:ex04_staged}, the vertical displacement at the tunnel crown increases progressively with the excavation step, and all three solutions show consistent deformation trends. The quantitative comparison in Tab.~\ref{tab:excavation_displacement_error} further indicates that the CS-FEM results are closer to the reference solution than the FEM results throughout the five excavation steps. The average relative errors calculated from the tabulated data are $1.28\times10^{-2}$ for FEM and $3.19\times10^{-3}$ for CS-FEM, demonstrating the improved accuracy of the proposed method in the staged excavation analysis. As shown in Fig.~\ref{fig:ex04_contour}, the vertical displacement contours obtained by CS-FEM and FEM agree well with the reference solution, indicating the feasibility of the proposed method for staged tunnel excavation analysis.

\begin{table}[H]
\centering
\caption{Mechanical parameters of the surrounding rock.}
\label{tab:ex04_rock_parameters}
\begin{tabular}{ccccc}
\hline
Parameter 
& $E$
& $\nu$
& $c$ 
& $\varphi$ \\
\hline
Value 
& $8~\mathrm{MPa}$ 
& $0.27$ 
& $110~\mathrm{kPa}$ 
& $39^\circ$ \\
\hline
\end{tabular}
\end{table}

\begin{table}[H]
\centering
\caption{Comparison of vertical displacements and relative errors at different excavation steps.}
\label{tab:excavation_displacement_error}
\resizebox{\textwidth}{!}{\begin{tabular}{cccccc}
\toprule
Excavation step & FEM & CS-FEM & Reference solution & FEM error & CS-FEM error \\
\midrule
1 & -0.12933 & -0.13331 & -0.13266 & $2.51\times10^{-2}$ & $4.90\times10^{-3}$ \\
2 & -0.15247 & -0.15366 & -0.15450 & $1.31\times10^{-2}$ & $5.44\times10^{-3}$ \\
3 & -0.15866 & -0.15976 & -0.16015 & $9.30\times10^{-3}$ & $2.44\times10^{-3}$ \\
4 & -0.16265 & -0.16376 & -0.16405 & $8.53\times10^{-3}$ & $1.77\times10^{-3}$ \\
5 & -0.16626 & -0.16738 & -0.16762 & $8.11\times10^{-3}$ & $1.43\times10^{-3}$ \\
\bottomrule
\end{tabular}}
\end{table}

\subsection{Slope stability analysis}
\subsubsection{Benchmark example}
\begin{figure}[h]
  \centering
  \includegraphics[width=1.0\textwidth]{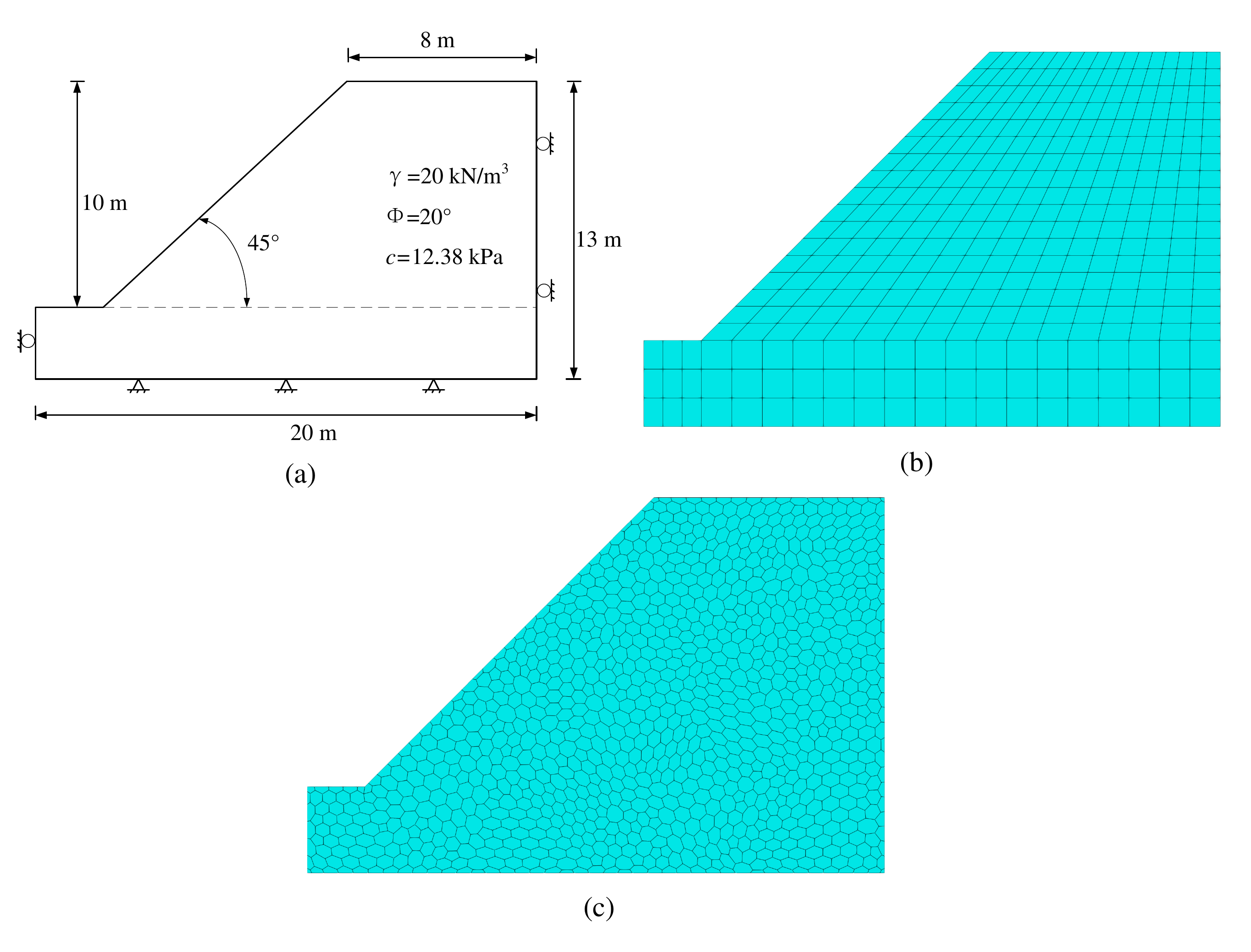}
  \caption{\textcolor{black}{Schematic diagram for the slope stability analysis; (a) geometry and boundary
conditions; (b) FEM mesh; (c) polygonal mesh}}
  \label{fig:slope}
\end{figure}

To verify the accuracy of the proposed method, we examined a uniform slope with a height of $H = 10 \text{m}$, sloping at an angle of $\beta = 45^\circ$ with a friction angle of $\phi = 20^\circ$, a unit weight of $\gamma = 20 \text{kN/m}^3$, a cohesion of $c = 12.38 \text{kPa}$, a Young's modulus of $E = 100 \text{MPa}$, and a Poisson's ratio of $\nu = 0.35$, as shown in Fig.~\ref{fig:slope} (a). Based on the soil properties, the slope has a factor of safety of exactly 1.0, as per the limit analysis solution presented by Chen~\citep{chen2012limit}. The lateral boundaries of the slope are constrained in the normal direction, while the bottom boundary is fixed in all directions.

The stability of slope is evaluated using the shear strength reduction technique \citep{matsui1992finite}. The parameters of shear strength are defined as 
\begin{equation}c_m=\frac c{F_r},\end{equation}
\begin{equation}\phi_m=\arctan\biggl(\frac{\tan\phi}{F_r}\biggr),\end{equation}
where $c$ and $\phi$ are the actual shear strength parameters, $c_m$ and $\phi_m$ are the shear strength which can maintain the stability of slope. $F_r$ is the strength reduction factor. In the reduction process, the input shear strength parameters $c$ and $\phi$ are reduced to $c_m$ and $\phi_m$ to trigger slope failure. The failure criterion was defined by the big jump in the nodal displacement at a chosen point close to the slope’s surface.

\begin{figure}[H]
  \centering
  \includegraphics[width=1.0\textwidth]{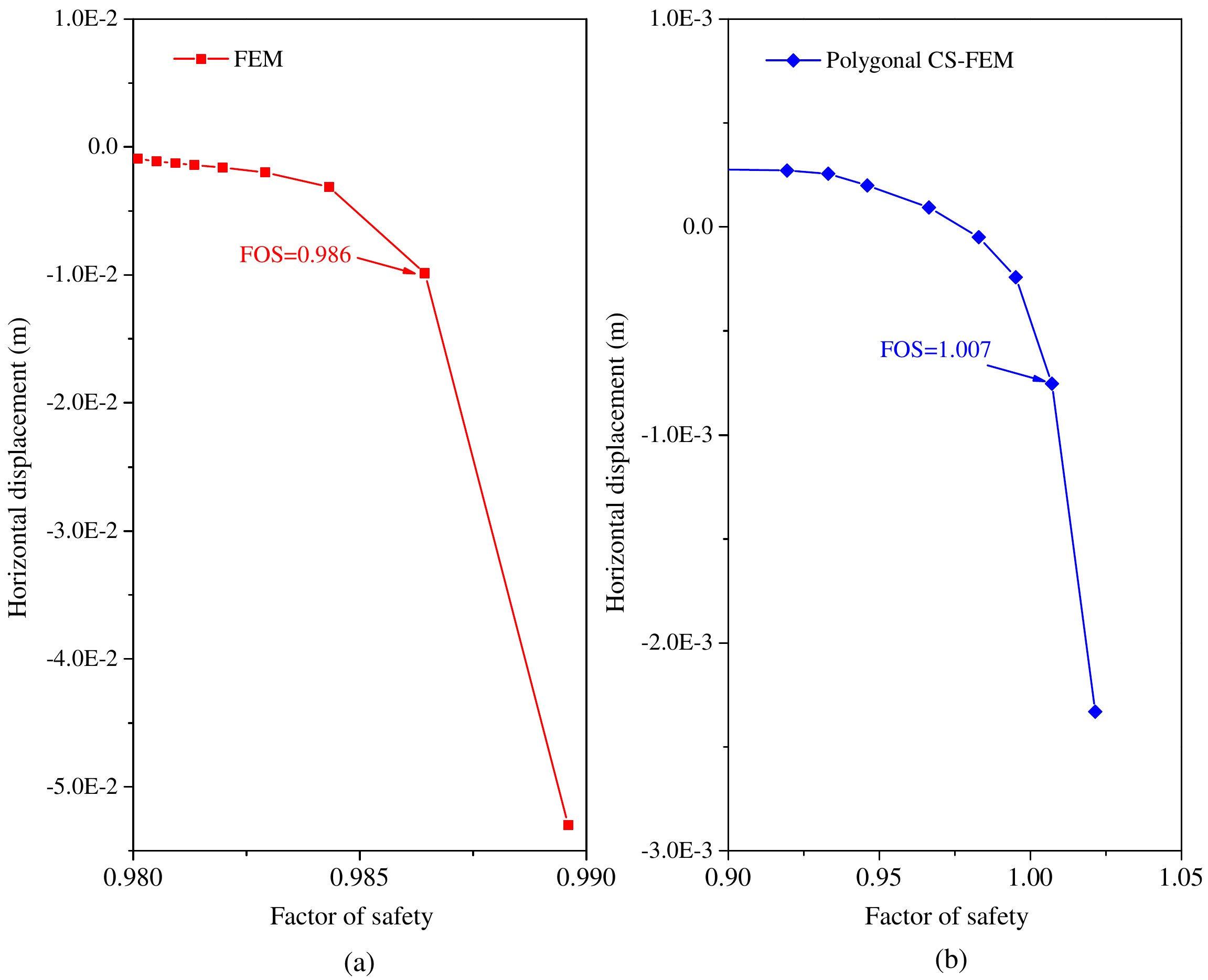}
  \caption{Relationship between horizontal displacement and the factor of safety; (a) FEM result; (b) polygonal CS-FEM result.}
  \label{fig:ex05_fos}
\end{figure}

\begin{figure}[H]
  \centering
  \includegraphics[width=1.0\textwidth]{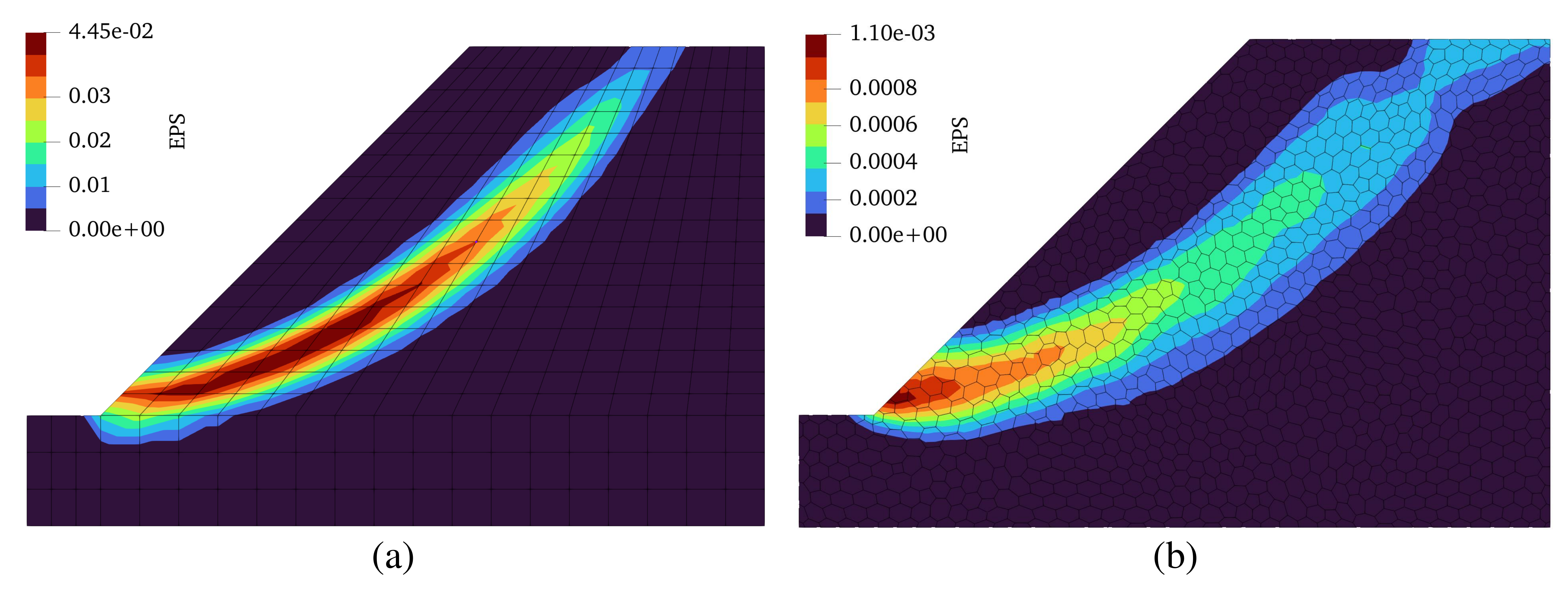}
  \caption{Plastic strain distribution in the slope at failure; (a) FEM result; (b) polygonal CS-FEM result.}
  \label{ex05:fig:slope_eps}
\end{figure}

We extracted the horizontal displacement of the slope vertex and plotted the relationship between horizontal displacement and the safety factor, as shown in Fig.~\ref{fig:ex05_fos}. The displacement–safety factor curves indicate that the critical safety factors obtained by CS-FEM and FEM are 1.007 and 0.986, respectively. Compared with the reference value of 1.0 reported in the literature~\citep{dawson1999slope}, the relative errors of CS-FEM and FEM are $7.0\times10^{-3}$ and $1.4\times10^{-2}$, respectively. This comparison shows that the polygonal CS-FEM provides a more accurate prediction of the slope safety factor. Moreover, Fig.~\ref{ex05:fig:slope_eps} shows the distribution of plastic strain in the slope at failure. It can be seen from the figure that the plastic zones have fully penetrated when the slope fails.

\subsubsection{Stability analysis of a slope with geological strata}

\begin{figure}[H]
  \centering
  \includegraphics[width=1.0\textwidth]{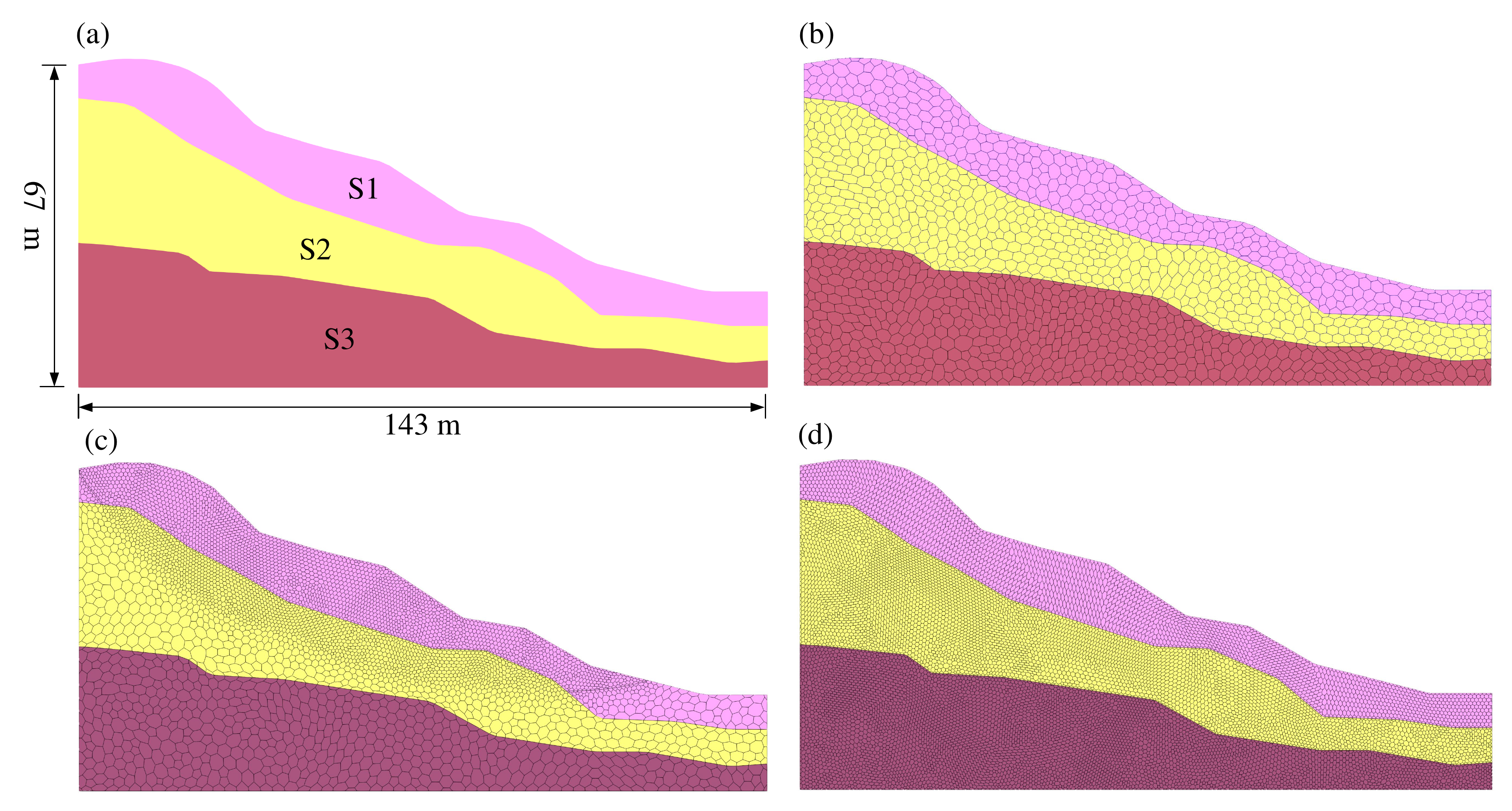}
  \caption{Geometric model and mesh discretizations for the slope case analysis: (a) slope geometry with three soil layers; (b) coarse mesh; (c) locally refined mesh; (d) fine mesh.}
  \label{ex05:slope_eng_geo_mesh}
\end{figure}

A layered slope is considered in this section to further evaluate the applicability of the proposed polygonal CS-FEM to practical slope stability analysis. The geometric model and three polygonal mesh discretizations are shown in Fig.~\ref{ex05:slope_eng_geo_mesh}. As shown in Fig.~\ref{ex05:slope_eng_geo_mesh}(a), the slope has a horizontal length of 143~m and a maximum height of 67~m. The computational domain consists of three soil layers, denoted as S1, S2, and S3, forming an irregular stratigraphic structure along the slope profile. To examine the mesh applicability and convergence behavior of the proposed method, three polygonal meshes are employed: a coarse mesh, a locally refined mesh, and a fine mesh, as shown in Fig.~\ref{ex05:slope_eng_geo_mesh}(b)--(d), respectively. The locally refined mesh is refined around the potential slip surface and the stratigraphic interfaces, where large gradients of displacement and plastic deformation are expected. These polygonal meshes can flexibly conform to the irregular layer interfaces and slope boundary, thereby avoiding excessive mesh distortion near the stratigraphic interfaces. The material parameters of the three strata are listed in Tab.~\ref{tab:slope_strata_parameters}.

\begin{figure}[H]
  \centering
  \includegraphics[width=1.0\textwidth]{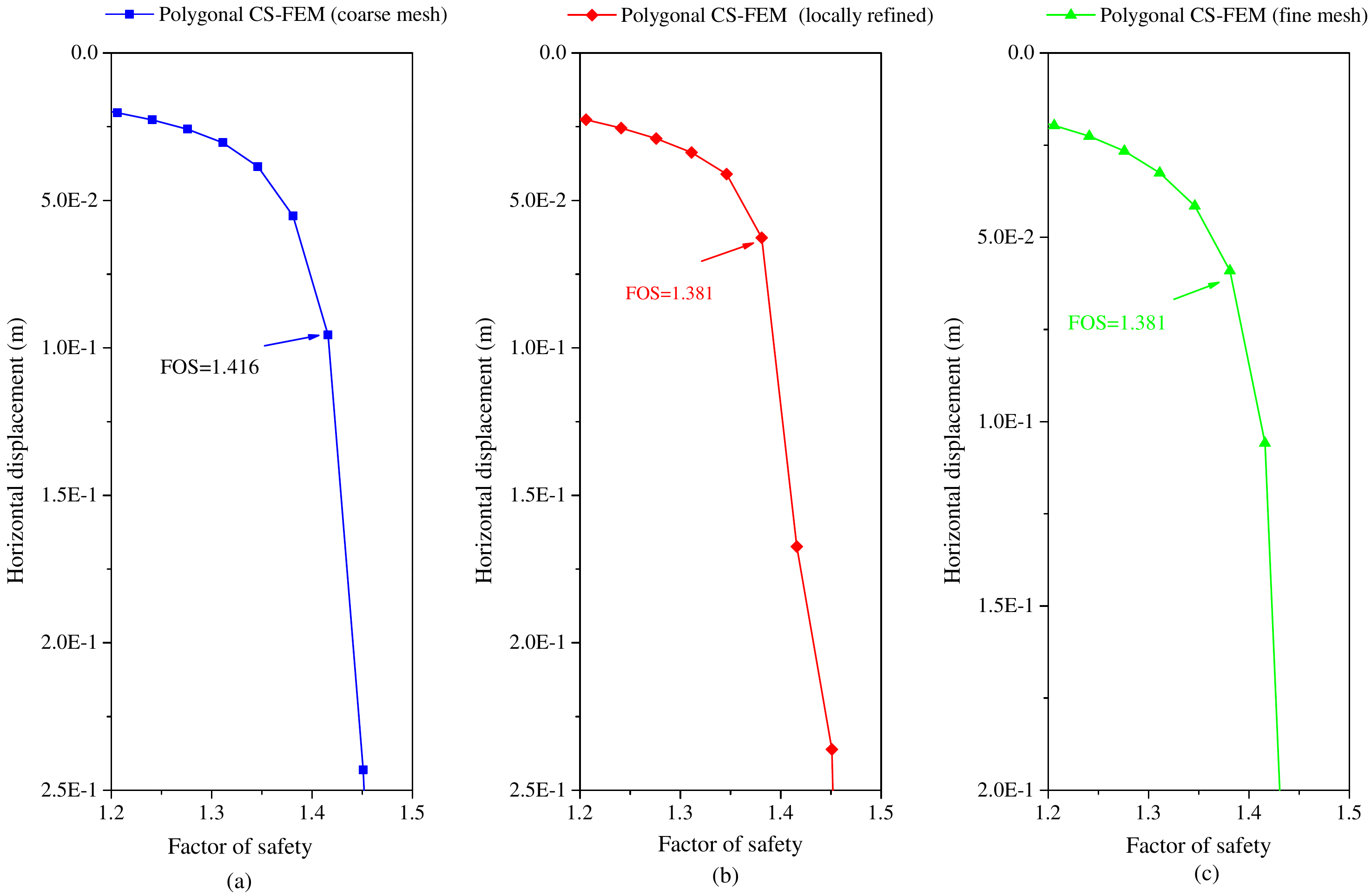}
\caption{Relationship between horizontal displacement and the factor of safety; (a) the result of coarse mesh; (b) the result of locally refined mesh; (c) the result of fine mesh.}
  \label{ex05:slope_eng_fos}
\end{figure}

The relationship between the horizontal displacement and the strength reduction factor is shown in Fig.~\ref{ex05:slope_eng_fos}. The factor of safety obtained using the coarse mesh is 1.416, whereas both the locally refined and fine meshes give a value of 1.381. This indicates that the locally refined mesh can achieve nearly the same stability result as the fine mesh.

\begin{figure}[H]
  \centering
  \includegraphics[width=0.65\textwidth]{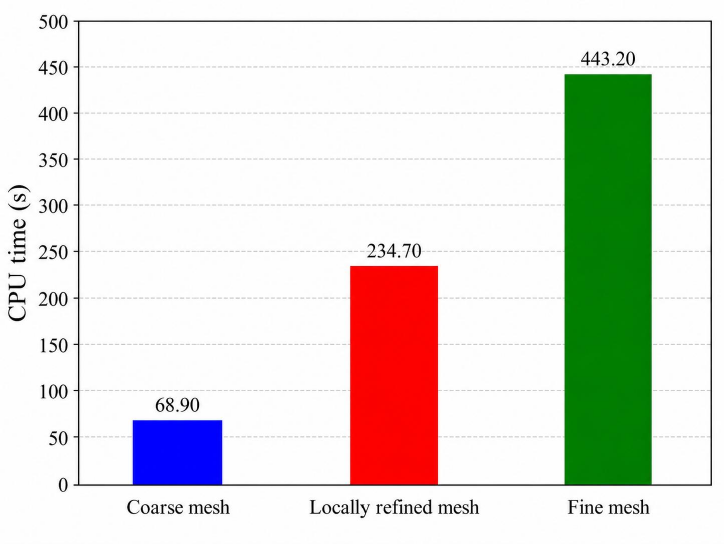}
  \caption{Comparison of CPU time for the slope stability analysis using different polygonal mesh discretizations.}
  \label{ex05:slope_eng_cpu_time}
\end{figure}

The computational efficiency of the three meshes is compared in Fig.~\ref{ex05:slope_eng_cpu_time}. The CPU time increases from 68.90~s for the coarse mesh to 234.70~s for the locally refined mesh and 443.20~s for the fine mesh. Compared with the fine mesh, the locally refined mesh reduces the CPU time by approximately 47.04\% while maintaining the same factor of safety, demonstrating a favorable balance between accuracy and computational efficiency.

\begin{figure}[H]
  \centering
  \includegraphics[width=1.0\textwidth]{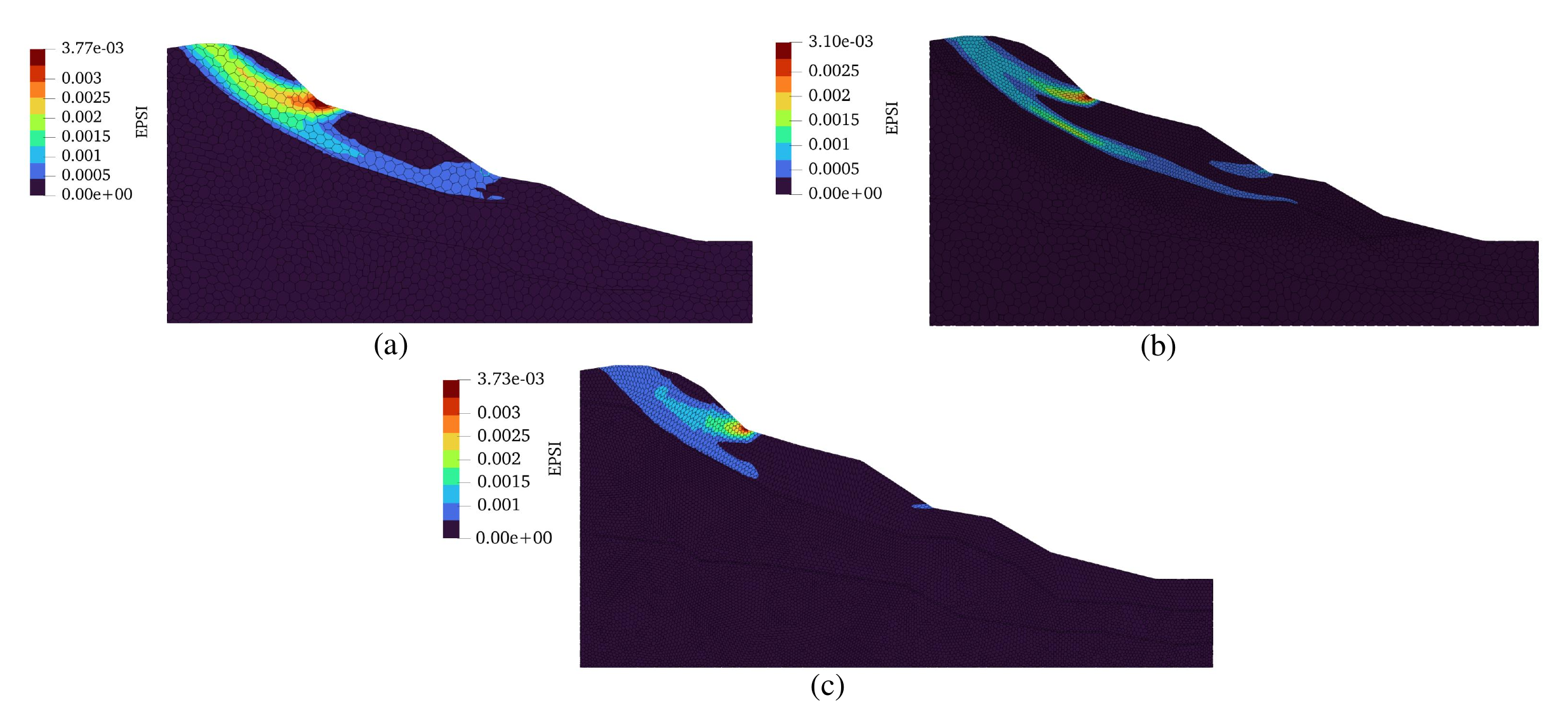}
\caption{Plastic strain distribution in the slope at failure; (a) coarse mesh; (b) locally refined mesh; (c) fine mesh.}
  \label{slope_eng_epi}
\end{figure}

Fig.~\ref{slope_eng_epi} presents the plastic strain distribution of the slope at failure. The plastic zone is mainly located in the upper part of the slope and along the potential slip surface, reflecting the failure mode of the layered slope. The locally refined mesh provides a plastic zone pattern close to that obtained from the fine mesh, further confirming its effectiveness for this slope stability analysis.

\begin{table}[H]
\centering
\caption{Material parameters of the three soil strata in the slope case analysis.}
\label{tab:slope_strata_parameters}
\begin{tabular}{cccccc}
\toprule
Stratum & $E$ & $\nu$ & $c$ & $\phi$ & $\gamma$ \\
        & (MPa) & (-) & (kPa) & ($^\circ$) & (kN/m$^3$) \\
\midrule
S1 & 98 & 0.32 & 210 & 31 & 21.3 \\
S2 & 322 & 0.29 & 563 & 42 & 21.8 \\
S3 & 546 & 0.27 & 649 & 43 & 22.3 \\
\bottomrule
\end{tabular}
\end{table}

\section{Conclusions}
\label{sec:conclusions}

This study developed a polygonal cell-based smoothed finite element method for nonlinear geotechnical analysis and implemented it in ABAQUS through the UEL interface. The proposed method combines Wachspress interpolation, cell-based strain smoothing, and incremental elasto-plastic constitutive updates. By evaluating the smoothed strain--displacement matrix through boundary integration, the formulation avoids direct calculation of shape-function derivatives inside polygonal elements and enables polygonal meshes and hybrid quadtree meshes with hanging nodes to be handled in a unified framework.

The accuracy and applicability of the proposed method were verified through several benchmark and engineering examples, including an infinite plate with a circular hole, strip footing, core rockfill dam, tunnel excavation, and slope stability problems. The numerical results show that the polygonal CS-FEM can accurately predict displacement, stress, plastic strain, bearing capacity, and factor of safety. Compared with conventional FEM, the proposed method provides better mesh flexibility and accurate nonlinear response prediction, especially for problems involving irregular geometries, local refinement, staged construction, and localized failure.

The locally refined polygonal CS-FEM further improves the balance between accuracy and efficiency. In the strip footing problem, it produced the smallest bearing-capacity error, and in the layered slope example, it obtained the same factor of safety as the fine mesh while reducing the CPU time by approximately 47.04\%. These results demonstrate that the proposed method is a reliable and efficient numerical tool for nonlinear geotechnical analysis. Future work will extend the formulation to three-dimensional polyhedral elements, coupled hydro-mechanical analysis, and more advanced constitutive models.

\section{Acknowledgements}
The Fundamental Research Funds for the Central Universities (grant NO. B240201147), the Yunnan Fundamental Research Projects (grant NO. 202401CF070043), and the Xing Dian Talent Support Program of Yunnan Province (XDYC-QNRC-2022-0764) provided support for this study.

 \bibliographystyle{elsarticle-harv} 
 \bibliography{cas-refs}




% %% For authoryear reference style

% \bibitem[Lamport(1994)]{lamport94}
%   \textit{\LaTeX: a document preparation system},
%   Addison Wesley, Massachusetts,

% \end{thebibliography}
\end{document}